# Adaptive Mollifiers – High Resolution Recovery of Piecewise Smooth Data from its Spectral Information

Eitan Tadmor     Jared Tanner

January 8, 2018

To Ron DeVore with Friendship and Appreciation

**Abstract**

We discuss the reconstruction of piecewise smooth data from its (pseudo-) spectral information. Spectral projections enjoy superior resolution provided the data is globally smooth, while the presence of jump discontinuities is responsible for spurious $\mathcal{O}(1)$ Gibbs oscillations in the neighborhood of edges and an overall deterioration to the unacceptable first-order convergence rate. The purpose is to regain the superior accuracy in the piecewise smooth case, and this is achieved by mollification.

Here we utilize a modified version of the two-parameter family of spectral mollifiers introduced by Gottlieb & Tadmor [GoTa85]. The ubiquitous one-parameter, finite-order mollifiers are based on *dilation*. In contrast, our mollifiers achieve their high resolution by an intricate process of high-order *cancelation*. To this end, we first implement a localization step using edge detection procedure, [GeTa00a, GeTa00b]. The accurate recovery of piecewise smooth data is then carried out in the direction of smoothness away from the edges, and *adaptivity* is responsible for the high resolution. The resulting adaptive mollifier greatly accelerates the convergence rate, recovering piecewise analytic data within exponential accuracy while removing spurious oscillations that remained in [GoTa85]. Thus, these adaptive mollifiers offer a robust, general-purpose "black box" procedure for accurate post processing of piecewise smooth data.

## Contents

# 1 Introduction

We study a new procedure for high resolution recovery of piecewise smooth data from its (pseudo-)spectral information. The purpose is to overcome the low-order accuracy and spurious oscillations associated with Gibbs phenomena, and to regain the superior accuracy encoded in the global spectral coefficients.

A standard approach for removing spurious oscillations is based on mollification over a *local* region of smoothness. To this end one employs a one-parameter family of dilated unit mass mollifiers of form $\varphi_\theta = \varphi(x/\theta)/\theta$. In general, such compactly supported mollifiers are restricted to finite-order accuracy, $|\varphi_\theta \star f(x) - f(x)| \leq C_r \theta^r$, depending on the number $r$ of vanishing moments $\varphi$ has. Convergence is guaranteed by letting the *dilation* parameter, $\theta \downarrow 0$.

In [GoTa85] we introduced a two-parameter family of spectral mollifiers of the form

$$\psi_{p,\theta}(x) = \frac{1}{\theta}\rho(\frac{x}{\theta})D_p(\frac{x}{\theta}).$$

Here $\rho(\cdot)$ is an arbitrary $C_0^\infty(-\pi,\pi)$ function which localizes the $p$-degree Dirichlet kernel $D_p(y) := \frac{\sin(p+1/2)y}{2\pi \sin(y/2)}$. The first parameter — the dilation parameter $\theta$ need not be small in this case, in fact $\theta = \theta(x)$ is made as large as possible while maintaining the smoothness of $\rho(x - \theta \cdot)f(\cdot)$. Instead, it is the second parameter – the degree $p$, which allows the high accuracy recovery of piecewise smooth data from its (pseudo-)spectral projection, $P_N f(x)$. The high accuracy recovery is achieved here by choosing large $p$'s, enforcing an intricate process of *cancelation* as an alternative to the usual finite-order accurate process of localization.

In §2 we begin by revisiting the convergence analysis of [GoTa85]. Spectral accuracy is achieved by choosing an increasing $p \sim \sqrt{N}$, so that $\psi_{p,\theta}$ has *essentially* vanishing moments all orders, $\int y^s \psi_{p,\theta}(y)dy = \delta_{s0} + C_s \cdot N^{-s/2}$, $\forall s$, yielding the 'infinite-order' accuracy bound in the sense of $|\psi_{p,\theta} \star P_N f(x) - f(x)| \leq C_s \cdot N^{-s/2}$, $\forall s$.

Although the last estimate yields the desired spectral convergence rate sought for in [GoTa85], it suffers as an over-pessimistic restriction since its derivation ignores the possible dependence of $p$ on the degree of local smoothness, $s$, and the support of local smoothness, $\sim \theta = \theta(x)$. In §3 we begin a detailed study on the optimal choice of the $(p,\theta)$ parameters of the spectral mollifiers, $\psi_{p,\theta}$:

- Letting $d(x)$ denote the distance to the nearest edge, we first set $\theta = \theta(x) \sim d(x)$ so that $\psi_{p,\theta} \star P_N f(x)$ incorporates the largest smooth neighborhood around $x$. To find the distance to the nearest discontinuity we utilize a general edge detection procedure, [GeTa99, GeTa00a, GeTa00b], where the location (and amplitudes) of all edges are found in *one global sweep*. Once the edges are located it is a straightforward matter to evaluate, at every $x$, the appropriate spectral parameter, $\theta(x) = d(x)/\pi$.

- Next, we turn to examine the degree $p$, which is responsible for the overall high accuracy by enforcing an intricate cancelation. A careful analysis carried out in §3.1 leads to an optimal choice of an *adaptive* degree of order $p = p(x) \sim d(x)N$. Indeed, numerical experiments reported back in the original [GoTa85] and additional tests carried out in §3.2 below and which motivate the present study, clearly indicate a superior convergence up to the immediate vicinity of the interior edges with an *adaptive* degree of the optimal order $p = p(x) \sim d(x)N$.

Given the spectral projection of a piecewise analytic function, $S_N f(\cdot)$, our 2-parameter family of adaptive mollifiers, equipped with the optimal parameterization outlined above yields – consult Theorem 3.1 below,

$$|\psi_{p,\theta} \star S_N f(x) - f(x)| \leq Const \cdot d(x) N \cdot e^{-\eta\sqrt{d(x)N}}.$$

The last error bound shows that the adaptive mollifier is exponentially accurate at all $x$'s except for the immediate $\mathcal{O}(1/N)$-neighborhood of the jumps of $f(\cdot)$ where $d(x) \sim 1/N$. We note in passing the rather remarkable dependence of this error estimate on the $C_0^\infty$ regularity of $\rho(\cdot)$. Specifically, the exponential convergence rate of a fractional power is related to the Gevrey regularity of the localizer $\rho(\cdot)$; in this paper we use the $G_2$-regular cut-off $\rho_c(x) = exp(cx^2/(x^2 - \pi^2))$ which led to the fractional power $1/2$.

Similar results holds in the discrete case. Indeed, in this case, one can bypass the discrete Fourier coefficients: expressed in terms of the given equidistant discrete values, $\{f(y_\nu)\}$, of piecewise analytic $f$, we have – consult Theorem 3.2 below,

$$|\frac{\pi}{N} \sum_{\nu=0}^{2N-1} \psi_{p,\theta}(x - y_\nu) f(y_\nu) - f(x)| \leq Const \cdot (d(x)N)^2 \cdot e^{-\eta\sqrt{d(x)N}}.$$

Thus, the discrete convolution $\sum_\nu \psi_{p,\theta}(x - y_\nu) f(y_\nu)$ forms an exponentially accurate near-by interpolant[1], which serves as an effective tool to reconstruct the intermediate values of piecewise smooth data. These near-by "expolants" are reminicient of quasi-interpolants, e.g., [BL93], with the emphasize given here to nonlinear adaptive recovery which is based on *global* regions of smoothness.

What happens in the immediate, $\mathcal{O}(1/N)$-neighborhood of the jumps? in §4 we complete our study of the adaptive mollifiers by introducing a novel procedure of *normalization*. Here we enforce the first few moments of the spectral mollifier, $\psi_p = \rho D_p$ to vanish, so that we regain *polynomial* accuracy in the immediate neighborhood of the jump. Taking advantage of the freedom in choosing the localizer $\rho(\cdot)$, we show how to modify $\rho$ to regain the local accuracy by enforcing finitely many vanishing moments of $\psi_p = \rho D_p$, while retaining the same overall exponential outside the immediate vicinity of the jumps. By appropriate normalization, the localized Dirichlet kernel we introduce maintains at least second order convergence *up to* the discontinuity. Increasingly higher orders of accuracy can be worked out as we move further away from these jumps and eventually turning into the exponentially accurate regime indicate earlier. In summary, the spectral mollifier amounts to a variable order recovery procedure adapted to the *number of cells* from the jump discontinuities, which is reminiscent of the variable order, Essentially Non Oscillatory piecewise polynomial reconstruction in [HEOC85]. The currrent procedure is also reminicient of the $h$-$p$ methods of Babuška and his collabrators, with the emphasize given here to an increasing number of *global* moments ($p$) without the "$h$"-refinement. The numerical experiments reported in §3.2 and §4.3 confirm the superior high resolution of the spectral mollifier $\psi_{p,\theta}$ equipped with the proposed optimal parameterization.

**Acknowledgment**. Research was supported in part by ONR Grant No. N00014-91-J-1076 and NSF grant #DMS01-07428.

## 2  Spectral Mollifiers

### 2.1  The two-parameter spectral mollifier $\psi_{p,\theta}$

The Fourier projection of a $2\pi$-periodic function $f(\cdot)$,

---
[1] Called expolant for short

$$S_N f(x) := \sum_{|k| \leq N} \hat{f}_k e^{ikx}, \qquad \hat{f}_k := \frac{1}{2\pi} \int_{-\pi}^{\pi} f(x) e^{-ikx} dx \qquad (2.1)$$

enjoys the well known spectral convergence rate, that is, the convergence rate is as rapid as the *global* smoothness of $f(\cdot)$ permits in the sense that for *any* $s$ we have[2]

$$|S_N f(x) - f(x)| \leq Const \|f\|_{C^s} \cdot \frac{1}{N^{s-1}} \qquad \forall s. \qquad (2.2)$$

Equivalently, this can be expressed in terms of the usual Dirichlet Kernel

$$D_N(x) := \frac{1}{2\pi} \sum_{k=-N}^{N} e^{ikx} \equiv \frac{\sin(N+1/2)x}{2\pi \sin(x/2)}, \qquad (2.3)$$

where $S_N f \equiv D_N \star f$, and the spectral convergence statement in (2.2) recast into the form

$$|D_N \star f(x) - f(x)| \leq Const \|f\|_{C^s} \cdot \frac{1}{N^{s-1}} \qquad \forall s. \qquad (2.4)$$

Furthermore, if $f(\cdot)$ is analytic with analyticity strip of width $2\eta$, then $S_N f(x)$ is characterized by an exponential convergence rate, e.g., [Ch, Ta94]

$$|S_N f(x) - f(x)| \leq Const_\eta \cdot N e^{-N\eta}. \qquad (2.5)$$

If, on the other hand, $f(\cdot)$ experiences a simple jump discontinuity, say at $x_0$, then $S_N f(x)$ suffers from the well known Gibbs' phenomena, where the uniform convergence of $S_N f(x)$ is lost in the neighborhood of $x_0$, and moreover, the *global* convergence rate of $S_N f(x)$ deteriorates to first order. To accelerate the slow convergence rate, we focus our attention on the classical process of mollification. Standard mollifiers are based on a one-parameter family of dilated unit mass functions of the form

$$\varphi_\theta(x) := \frac{1}{\theta} \varphi\left(\frac{x}{\theta}\right) \qquad (2.6)$$

which induce convergence by letting $\theta$ to zero. In general, $|\varphi_\theta \star f(x) - f(x)| \leq C_r \theta^r$ describes the convergence rate of *finite* order $r$, where $\varphi$ possesses $r$ vanishing moments

$$\int y^s \varphi(y) dy = \delta_{s0} \qquad s = 0, 1, 2, \ldots, r-1. \qquad (2.7)$$

In the present context of recovering *spectral* convergence, however, we follow Gottlieb and Tadmor, [GoTa85], using a two-parameter family of mollifiers, $\psi_{p,\theta}(x)$, where $\theta$ is a dilation parameter, $\psi_{p,\theta}(x) = \psi_p(x/\theta)/\theta$, and $p$ stipulates how closely $\psi_{p,\theta}(x)$ possesses near vanishing moments. To form $\psi_p(x)$, we let $\rho(x)$ be an arbitrary $C_0^\infty$ function supported in $(-\pi, \pi)$ and we consider the localized Dirichlet kernel

$$\psi_p(x) := \rho(x) D_p(x). \qquad (2.8)$$

Our two-parameter mollifier is then given by the dilated family of such localized Dirichlet kernels

$$\psi_{p,\theta}(x) := \frac{1}{\theta} \psi_p\left(\frac{x}{\theta}\right) \equiv \frac{1}{\theta} \rho\left(\frac{x}{\theta}\right) D_p\left(\frac{x}{\theta}\right). \qquad (2.9)$$

---

[2]Here and below we denote the usual $\|f\|_{C^s} := \|f^{(s)}\|_{L^\infty}$.

According to (2.9), $\psi_{p,\theta}$ consists of two ingredients, $\rho(x)$ and $D_p(x)$, each has essentially separate role associated with the two independent parameters $\theta$ and $p$. The role of $\rho\left(\frac{x}{\theta}\right)$ is, through its $\theta$-dependence, to *localize* the support of $\psi_{p,\theta}(x)$ to $(-\theta\pi, \theta\pi)$. The Dirichlet kernel $D_p(x)$ is charged, by varying $p$, with controlling the increasing number of near vanishing moments of $\psi_{p,\theta}$, and hence the overall superior accuracy of our mollifier. Indeed, by imposing the normalization of

$$\rho(0) = 1, \tag{2.10}$$

we find that an increasing number of moments of $\psi_{p,\theta}$ are of the vanishing order $\mathcal{O}(p^{-(s-1)})$,

$$\int_{-\pi\theta}^{\pi\theta} y^s \psi_{p,\theta}(y) dy = \int_{-\pi}^{\pi} (y\theta)^s \rho(y) D_p(y) dy = D_p \star (y\theta)^s \rho(y)_{|y=0} = \delta_{s0} + C_s \cdot p^{-(s-1)} \quad \forall s, \tag{2.11}$$

where according to (2.4), $C_s = Const\|(y\theta)^s\rho(y)\|_{C^s}$. We shall get into a detailed convergence analysis in the discussion below.
We conclude this section by highlighting the contrast between the standard, polynomially accurate mollifier, (2.7) and the spectral mollifiers (2.9). The former depends on one dilation parameter, $\theta$, which is charged of inducing a fixed order of accuracy by letting $\theta \downarrow 0$. Thus, in this case convergence is enforced by *localization*, which is inherently limited to a fixed polynomial order. The spectral mollifier, however, has the advantage of employing two free parameters: the dilation parameter $\theta$ which need not be small - in fact, $\theta$ is made as *large* as possible while maintaining $\rho(x - \theta y) f(y)$ free of discontinuities; the need for this desired smoothness will be made more evident in the next section. It is the second parameter, $p$, which is in charge of enforcing the high accuracy by letting $p \uparrow \infty$. Here, convergence is enforced by a delicate process of *cancelation* which will enable us to derive, in §3, exponential convergence.

## 2.2 Error analysis for spectral mollifier

We now turn to consider the error of our mollification procedure, $E(N, p, \theta; f(x))$, at an arbitrary fixed point $x \in [0, 2\pi)$

$$E(N, p, \theta; f(x)) = E(N, p, \theta) := \psi_{p,\theta} \star S_N f(x) - f(x), \tag{2.12}$$

where we highlight the dependence on three free parameters at our disposal – the degree of the projection, N, the support of our mollifier, $\theta$, and the degree with which we approximate an arbitrary number of vanishing moments, $p$. The dependence on the degree of piecewise smoothness of $f(\cdot)$ will play a secondary role in the choice of these parameters.
We begin by decomposing the error into the three terms

$$E(N, p, \theta) = (f \star \psi_{p,\theta} - f) + (S_N f - f) \star (\psi_{p,\theta} - S_N \psi_{p,\theta}) + (S_N f - f) \star S_N \psi_{p,\theta}. \tag{2.13}$$

The last term, $(S_N f - f) \star S_N \psi_{p,\theta}$, vanishes by orthogonality, and hence we are left with the first and second terms, which we refer to as the Regularization and Truncation errors, respectively

$$E(N, p, \theta) \equiv (f \star \psi_{p,\theta} - f) + (S_N f - f) \star (\psi_{p,\theta} - S_N \psi_{p,\theta}) =: R(N, p, \theta) + T(N, p, \theta). \tag{2.14}$$

Sharp error bounds for the regularization and truncation errors were originally derived in [GoTa85], and a short re-derivation now follows.
For the regularization error we consider the function

$$g_x(y) := f(x - \theta y)\rho(y) - f(x) \tag{2.15}$$

where $f(x)$ is the fixed point value to be recovered through mollification. Applying (2.4) to $g_x(\cdot)$ while noting that $g_x(0) = 0$, then the regularization error does not exceed

$$\begin{aligned} |R(N, p, \theta)| &:= |f \star \psi_{p,\theta} - f| = \left| \int_{-\pi}^{\pi} [f(x - \theta y)\rho(y) - f(x)] D_p(y) dy \right| \\ &= |D_p \star g_x(y)|_{y=0}| = |(S_p g_x(y) - g_x(y))_{|y=0}| \leq Const. \|g_x(y)\|_{C^s} \cdot \frac{1}{p^{s-1}}. \end{aligned} \tag{2.16}$$

Applying Leibnitz rule to $g_x(y)$,

$$|g_x^{(s)}(y)| \leq \sum_{k=0}^{s} \binom{s}{k} \theta^k |f^{(k)}(x - \theta y)| \cdot |\rho^{(s-k)}(y)| \leq \|\rho\|_{C^s} \|f^{(s)}\|_{L^\infty_{loc}} (1 + \theta)^s, \tag{2.17}$$

gives the desired upper bound

$$|R(N, p, \theta)| \leq Const. \|\rho\|_{C^s} \|f^{(s)}\|_{L^\infty_{loc}} \cdot p \left( \frac{2}{p} \right)^s. \tag{2.18}$$

Here and below $Const$ represents (possibly different) generic constants; also, $\|\cdot\|_{L^\infty_{loc}}$ indicates the $L^\infty$ norm to be taken over the *local* support of $\psi_{p,\theta}$. Note that $\|f^{(s)}\|_{L^\infty_{loc}} < \infty$ as long as $\theta$ is chosen so that $f(\cdot)$ is free of discontinuities in $(x - \theta\pi, x + \theta\pi)$.

To upperbound the truncation error we use Young's inequality followed by (2.4),

$$\begin{aligned} |T(N, p, \theta)| &\leq \|(S_N f - f) \star (\psi_{p,\theta} - S_N \psi_{p,\theta})\|_{L^\infty} \\ &\leq \|S_N f - f\|_{L^1} \cdot \|\psi_{p,\theta} - S_N \psi_{p,\theta}\|_{L^\infty} \leq M \|S_N f - f\|_{L^1} \cdot \|\psi_{p,\theta}\|_{C^s} \frac{1}{N^{s-1}}. \end{aligned} \tag{2.19}$$

Leibnitz rule yields,

$$|\psi_{p,\theta}^{(s)}| \leq \theta^{-(s+1)} \sum_{k=0}^{s} \binom{s}{k} |\rho^{(s-k)}| \cdot |D_p^{(k)}| \leq \|\rho\|_{C^s} \left( \frac{1+p}{\theta} \right)^{s+1}, \tag{2.20}$$

and together with (2.19) we arrive at the upper bound

$$|T(N, p, \theta; f)| \leq Const \|S_N f - f\|_{L^1} \cdot \|\rho\|_{C^s} \cdot \frac{(1+p)N}{\theta} \left( \frac{1+p}{N\theta} \right)^s. \tag{2.21}$$

A slightly tighter estimate is obtained by replacing the $L^1 - L^\infty$ bounds with $L^2$ bounds for $f$'s with bounded variation,

$$\begin{aligned} |T(N, p, \theta)| &\leq \|S_N f - f\|_{L^2} \times \|S_N \psi_{p,\theta} - \psi_{p,\theta}\|_{L^2} \leq \\ &\leq Const. \|f\|_{BV} \cdot N^{-1/2} \times \|\psi_{p,\theta}^{(s)}\|_{L^2} \cdot N^{-(s-1/2)}, \end{aligned} \tag{2.22}$$

and (2.20) then yields

$$|T(N, p, \theta)| \leq Const. \|\rho\|_{C^s} \cdot N \left( \frac{1+p}{N\theta} \right)^{s+1}. \tag{2.23}$$

Using this together with (2.18), we conclude with an error bound of $E(N, p, \theta; f(x))$,

$$|\psi_{p,\theta} \star S_N f(x) - f(x)| \leq Const \|\rho\|_{C^s} \left[ N \left( \frac{1+p}{N\theta} \right)^{s+1} + p \left( \frac{2}{p} \right)^s \|f^{(s)}\|_{L^\infty_{loc}(x)} \right], \quad \forall s, \quad (2.24)$$

where $\|f^{(s)}\|_{L^\infty_{loc}} = \sup_{y \in (x-\theta\pi, x+\theta\pi)} |f^{(s)}|$ measures the local regularity of $f$. It should be noted that one can use different orders of degrees of smoothness, say an $r$ order of smoothness for the truncation and $s$ order of smoothness for the regularization, yielding

$$|E(N, p, \theta; f(x))| \leq Const. \left[ \|\rho\|_{C^r} \cdot N \left( \frac{1+p}{N\theta} \right)^{r+1} + \|\rho\|_{C^s} \cdot p \left( \frac{2}{p} \right)^s \|f^{(s)}\|_{L^\infty_{loc}} \right], \quad \forall r, s. \quad (2.25)$$

### 2.3 Fourier interpolant - error analysis for pseudospectral mollifier

The Fourier interpolant of a $2\pi$-periodic function, $f(\cdot)$, is given by

$$I_N f(y) := \sum_{|k| \leq N} \tilde{f}_k e^{iky}, \quad \tilde{f}_k := \frac{1}{2N} \sum_{\nu=0}^{2N-1} f(y_\nu) e^{-iky_\nu}. \quad (2.26)$$

We observe that the moments computed in the spectral projection (2.1) are replaced here by the corresponding trapezoidal rule evaluated at the equidistant nodes $y_\nu = \frac{\pi}{N}\nu$, $\nu = 0, 1, \ldots, 2N-1$. It should be noted that this approximation by the trapezoidal rule converts the Fourier-Galerkin projection to a Pseudo Spectral Fourier collocation (interpolation) representation. It is well known that the Fourier Interpolant also enjoys spectral convergence, i.e.

$$|I_N f(x) - f(x)| \leq Const \|f\|_{C^s} \cdot \frac{1}{N^{s-1}}, \quad \forall s. \quad (2.27)$$

Furthermore, if $f(\cdot)$ is analytic with analyticity strip of width $2\eta$, then $S_N f(x)$ is characterized by an exponential convergence rate [Ta94]

$$|S_N f(x) - f(x)| \leq Const_\eta \cdot N e^{-N\eta}. \quad (2.28)$$

If, however, $f(\cdot)$ experiences a simple jump discontinuity, then the Fourier Interpolant suffers from the reduced convergence rate similar to the Fourier projection. To accelerate the slowed convergence rate we again make use of our *two-parameter* mollifier (2.9). When convolving $I_N f(x)$ by our two parameter mollifier we approximate the convolution by the Trapezoidal summation

$$\psi_{p,\theta} \star I_N f(x) \sim \frac{\pi}{N} \sum_{\nu=0}^{2N-1} f(y_\nu) \psi_{p,\theta}(x - y_\nu). \quad (2.29)$$

We note that the summation in (2.29) bypasses the need to compute the pseudo spectral coefficients $\tilde{f}_k$. Thus, in contrast to the spectral mollifiers carried out in the Fourier Space [MMO78], we are able to work directly in the physical space through using the sampling of $f(\cdot)$ at the equidistant points, $f(y_\nu)$.

The resulting error of our discrete mollification at the fixed point $x$ is given by

$$E(N,p,\theta) := \frac{\pi}{N} \sum_{\nu=0}^{2N-1} f(y_\nu)\psi_{p,\theta}(x - y_\nu) - f(x). \tag{2.30}$$

As before, we decompose the error into two components

$$\begin{aligned} E(N,p,\theta) &= \left(\frac{\pi}{N} \sum_{\nu=0}^{2N-1} f(y_\nu)\psi_{p,\theta}(x - y_\nu) - f \star \psi_{p,\theta}\right) + (f \star \psi_{p,\theta} - f) \\ &=: A(N,p,\theta) + R(N,p,\theta), \end{aligned} \tag{2.31}$$

where $R(N,p,\theta)$ is the familiar regularization error, and $A(N,p,\theta)$ is the so-called aliasing error committed by approximating the convolution integral by a Trapezoidal sum. It can be shown that, for any $m > 1/2$, the aliasing error does not exceed the truncation error, e.g., [Ta94, (2.2.16)],

$$\|A(N,p,\theta)\|_{L^\infty} \leq M_m \|T(N,p,\theta; f^{(m)})\|_{L^\infty} \cdot N^{(1/2-m)}, \qquad m > 1/2. \tag{2.32}$$

We choose $m = 1$: inserting this into (2.21) with $f$ replaced by $f'$, and noting that $\|S_N f' - f'\|_{L^1} \leq Const \|f\|_{BV}\sqrt{N}$, we recover the same truncation error bound we had in (2.21),

$$\|A(N,p,\theta)\|_{L^\infty} \leq Const.|T(N,p,\theta; f')|\frac{1}{\sqrt{N}} \leq Const.\|\rho\|_{C^s} \cdot N^2 \left(\frac{1+p}{N\theta}\right)^{s+1}. \tag{2.33}$$

Consequently, the error after discrete mollification of the Fourier Interpolant satisfies the same bound as the mollified Fourier projection

$$|E(N,p,\theta; f(x))| \leq Const \|\rho\|_{C^s} \left[N^2 \left(\frac{1+p}{N\theta}\right)^{s+1} + p\left(\frac{2}{p}\right)^s \|f^{(s)}\|_{L^\infty_{loc}}\right], \qquad \forall s \geq 1/2. \tag{2.34}$$

We close by noting that the spectral and pseudospectral error bounds, (2.24) and (2.34), are of the exact same order. And as before, one can use different orders of degrees of smoothness for the regularization and aliasing errors.

## 2.4    On the choice of the $(\theta, p)$ parameters – spectral accuracy.

We now turn to asses the role of the parameters, $\theta$ and $p$, based on the spectral and pseudospectral error bounds (2.24) and (2.34). We first address the localization parameter $\theta$. According to the first term on the right of (2.24), and respectively – (2.34), the truncation, and respectively – aliasing error bounds decrease for increasing $\theta$'s. Thus we are motivated to choose $\theta$ as large as possible. However, the silent dependence on $\theta$ of the regularization error term in (2.24) and (2.34) appears through the requirement of localized regularity, i.e. $\|f^{(s)}\|_{L^\infty_{loc}} = \sup_{y \in (x-\theta\pi, x+\theta\pi)} |f^{(s)}(y)| < \infty$. Hence, if $d(x)$ denotes the distance from $x$ to the nearest jump discontinuity of $f$,

$$d(x) := dist(x, sing\ supp\ f), \tag{2.35}$$

we then set

$$\theta := \frac{d(x)}{\pi} \leq 1. \tag{2.36}$$

This choice of $\theta$ provides us with the largest admissible support of the mollifier $\psi_{p,\theta}$, so that $\psi_{p,\theta} * f(x)$ incorporates only the (largest) smooth neighborhood around $x$. This results in an *adaptive* mollifier which amounts to a symmetric windowed filter of maximal width, $2d(x)$, to be carried out in the physical space. We highlight the fact that this choice of an $x$-dependent, $\theta(x) = d(x)/\pi$, results in a spectral mollifier that is *not* translation invariant. Consequently, utilizing such an adaptive mollifier is quite natural in the physical space, and although possible, it is not well suited for convolution in the frequency space.

How can we find the nearest discontinuity? we refer the reader to [GeTa99, GeTa00a, GeTa00b], where a general procedure to detect the edges in piecewise smooth data from its (pseudo-)spectral content. The procedure – carried in the physical space, is based on appropriate choice of concentration factors which lead to (generalized) conjugate sums which tend to concentrate in the vicinity of edges and are vanishing elsewhere. The locations (and amplitudes) of all the discontinuous jumps are found in *one global sweep*. Equipped with these locations, it is a straightforward matter to evaluate, at every $x$, the appropriate spectral parameter, $\theta(x) = d(x)/\pi$.

Next we address the all important choice of $p$ which controls how closely $\psi_{p,\theta}$ possesses near vanishing moments of increasing order, (2.11). Before determining an optimal choice of $p$ let us revisit the original approach taken by Gottlieb and Tadmor [GoTa85]. To this end, we first fix an arbitrary degree of smoothness $s$, and focus our attention on the optimal dependence of $p$ *solely* on $N$. With this in mind, the dominant terms of the error bounds (2.24) and (2.34), are of order $(p/N)^s$ and $p^{-s}$, respectively. Equilibrating these competing terms gives $p = \sqrt{N}$, which results in the spectral convergence rate sought for in [GoTa85], namely, for an arbitrary $s$

$$|E(N,p,\theta)|_{p=\sqrt{N}} \leq Const_{s,\theta} N^{-s/2}. \tag{2.37}$$

Although this estimate yields the desired spectral convergence rate sought for in [GoTa85], it suffers as an over-pessimistic restriction since the possible dependence of $p$ on $s$ and $\theta$ were not fully exploited. In fact, while the above approach of equilibration with $p$-depending solely on $N$ yields $p = N^{0.5}$, numerical experiments reported back in the original [GoTa85] have shown that when treating $p$ as a fixed power of $N$, $p = N^\beta$, superior results are obtained for $0.7 < \beta < 0.9$. Indeed, the numerical experiments reported in §3.2 below and which motivate the present study, clearly indicate that the contributions of the truncation and regularization terms are equilibrated when $p \sim N$. Moreover, the truncation and aliasing error contributions to the error bounds (2.24) and (2.34) predict convergence only for $x$'s which are bounded away from the jump discontinuities of $f$, where $\theta(x) > p/N$. Consequently, with $\theta(x) := d(x)/\pi$ and $f(\cdot)$ having a discontinuity, say at $x_0$, convergence can *not* be guaranteed in the region

$$(x_0 - \frac{p}{N}\pi, x_0 + \frac{p}{N}\pi). \tag{2.38}$$

Thus, a non-adaptive choice of $p$ – chosen as a fixed *fractional* power of $N$ *independent of* $\theta(x)$, say $p \sim \sqrt{N}$, can lead to a loss of convergence in a large zones of size $\mathcal{O}(N^{-1/2})$, around the discontinuity. The loss of convergence was confirmed in the numerical experiments reported in S3.2. This should be contrasted with the adaptive mollifiers introduced in the next §3, which will enable us to achieve exponential accuracy up to the immediate, $\mathcal{O}(1/N)$ vicinity of these discontinuities. We now turn to determine an optimal choice of $p$ by incorporating both – the distance to the nearest discontinuity, $d(x)$, and by exploiting the fact that the error bounds (2.24) and (2.34) allow us to use a variable degree of smoothness, $s$.

## 3  Adaptive Mollifiers – Exponential Accuracy

**Epilogue - Gevrey regularity**. The spectral decay estimates (2.2) and (2.27) tell us that for $C_0^\infty$ data, the (pseudo) spectral errors decay faster than any fixed polynomial order. To quantify the *actual* error decay, we need to classify the specific order of $C_0^\infty$ regularity. The Gevrey class, $G_\alpha, \alpha \geq 1$, consists of $\rho$'s with constants $\eta := \eta_\rho$ and $M := M_\rho$, such that the following estimate holds,

$$\sup_{x \in \Re} |\rho^{(s)}(x)| \leq M \frac{(s!)^\alpha}{\eta^s}, \qquad s = 1, 2, \ldots \tag{3.1}$$

We have two prototypical examples in mind.

**Example 1** A bounded analytic function $\rho$ belongs to $G_1$ with $M_\rho = \sup_{x \in \Re} |\rho(x)|$ and $2\eta_\rho$ equals the width of $\rho$'s analyticity strip.

**Example 2** Consider a $C_0^\infty(-\pi, \pi)$ cut-off function depending on an arbitrary constant $c > 0$, which takes the form

$$\rho_c(x) = \begin{cases} e^{\left(\frac{cx^2}{x^2 - \pi^2}\right)} & |x| < \pi \\ 0 & |x| \geq \pi \end{cases} =: e^{\left(\frac{cx^2}{x^2 - \pi^2}\right)} 1_{[-\pi, \pi]} \tag{3.2}$$

In this particular case there exists a constant $\lambda = \lambda_c$ such that the higher derivatives are upper bounded by[3]

$$|\rho_c^s(x)| \leq M \frac{s!}{(\lambda_c |x^2 - \pi^2|)^s} e^{\left(\frac{cx^2}{x^2 - \pi^2}\right)}, \qquad s = 1, 2, \ldots \tag{3.3}$$

The maximal value of the upper bound on the right hand side of (3.3) is obtained at $x = x_{max}$ where $x_{max}^2 - \pi^2 \sim -\pi^2 c/s$; consult[4]. This implies that our cut-off function $\rho_c$ admits $G_2$ regularity, namely, there exists a constant $\eta_c := \lambda_c \pi^2 c$ such that

$$\sup_{x \in \Re} |\rho_c^{(s)}(x)| \leq Const_c \cdot s! \left(\frac{s}{\eta_c}\right)^s e^{-s} \leq Const_c \frac{(s!)^2}{\eta_c^s} \qquad s = 1, 2, \ldots \tag{3.4}$$

We now turn to examine the actual decay rate of Fourier projections, $|S_N \rho - \rho|$, for arbitrary $G_\alpha$-functions. According to (2.2) combined with the growth of $\|\rho\|_{C^s}$ dictated by (3.1), the $L^\infty$ error in spectral projection of a $G_\alpha$ function, $\rho$, is governed by

$$|S_N \rho(x) - \rho(x)| \leq Const.N \frac{(s!)^\alpha}{(\eta N)^s}, \qquad s = 1, 2, \ldots \tag{3.5}$$

The expression of the type encountered on the right of (3.5), $(s!)^\alpha (\eta N)^{-s}$, attains its minimum at $s_{min} = (\eta N)^{1/\alpha}$,

---

[3] To this end note that $\rho_c(x) = e_+(x)e_-(x)$ with $e_\pm(x) := exp(cx/(x \pm \pi))$ for $x \in (-\pi, \pi)$. The functions $e_\pm(x)$ upper bounded by $|e_\pm^{(s)}(x)| \leq M_\pm s! (\lambda_c |x \pm \pi|)^{-s} e_\pm(x)$, with appropriate $\lambda = \lambda_c$, [Jo, p. 73]

[4] For large values of $s$, the function $|a(x)|^{-s} \cdot exp(\alpha a(x) + \beta/a(x))$ with fixed $\alpha$ and $\beta$ is maximized at $x = x_{max}$ such that $a(x_{max}) \sim -\beta/s$. In our case, $a(x) = x^2 - \pi^2$ and $\beta \sim c\pi^2$

$$\min_s \frac{(s!)^\alpha}{(\eta N)^s} \sim \min_s \left(\frac{s^\alpha}{\eta e^\alpha N}\right)^s = e^{-\alpha(\eta N)^{1/\alpha}}. \tag{3.6}$$

Thus, minimizing the upper-bound in (3.5) at $s = s_{min} = (\eta N)^{1/\alpha}$, yields the exponential accuracy of *fractional* order

$$|S_N \rho(x) - \rho(x)| \leq Const \cdot N e^{-\alpha(\eta N)^{1/\alpha}}, \qquad \rho \in G_\alpha. \tag{3.7}$$

The case $\alpha = 1$ recovers the exponential decay for analytic $\rho$'s, (2.5), whereas for $\alpha > 1$ we have exponential decay of fractional order. For example, our $G_2$ cut-off function $\rho = \rho_c$ in (3.2) satisfies (3.7) with $(\eta, \alpha) = (\eta_c, 2)$, yielding

$$|S_p \rho_c(x) - \rho_c(x)| \leq Const.p \cdot e^{-2\sqrt{\eta_c p}}. \tag{3.8}$$

Equipped with these estimates we now revisit the error decay of spectral mollifiers based on $G_\alpha$ cut-off functions $\rho$. Both contributions to the error in (2.14) — the regularization $R(N, p, \theta)$, and the truncation $T(N, p, \theta)$ (as well as aliasing $A(N, p, \theta)$ in (2.31)), are controlled by the decay rate of Fourier projections.

### 3.1 The $(\theta, p)$ parameters revisited – exponential accuracy.

We assume that $f(\cdot)$ is piecewise analytic. For each fixed $x$, our choice of $\theta = \theta(x) = d(x)/\pi$ guarantees that $f(x - \theta y)$ is analytic in the range $|y| \leq \pi$ and hence its product with the $G_\alpha(-\pi, \pi)$ function $\rho(y)$ yields the $G_\alpha$ regularity of $g_x(y) = f(x - \theta y)\rho(y) - f(x)$. According to (2.16), the regularization error, $R(N, p, \theta)$ is controlled by the Fourier projection of $g_x(\cdot)$, and in view of its $G_\alpha$ regularity, (3.7) yields

$$|R(N, p, \theta)| = |(S_p g_x(y) - g_x(y))_{|y=0}| \leq Const_\rho \cdot p \, e^{-\alpha(\eta p)^{1/\alpha}}. \tag{3.9}$$

For example, if $\rho = \rho_c$ we get

$$|R_{\rho_c}(N, p, \theta)| \leq Const_c \cdot p \cdot e^{-2\sqrt{\eta_c p}}. \tag{3.10}$$

*Remark.* It is here that we use the normalization, $\rho(0) - 1 = g_x(y = 0) = 0$, and (3.9) shows that one can slightly relax this normalization within the specified error bound

$$|\rho(0) - 1| \leq Const. \, e^{-\alpha(\eta p)^{1/\alpha}}. \tag{3.11}$$

Next we turn to the truncation error, $T(N, p, \theta)$. According to (2.19), its decay is controlled by the Fourier projection of the localized Dirichlet kernel $\psi_{p,\theta}(x) = \frac{1}{\theta}\psi_p\left(\frac{x}{\theta}\right)$. Here we shall need the specific structure of the localizer $\rho(x) = \rho_c(x)$ in (3.2). Leibnitz rule and (3.3) yield

$$\begin{aligned} |\psi_p^{(s)}(x)| &\leq \sum_{k=0}^s \binom{s}{k} |\rho_c^{(k)}(x)| \cdot |D_p^{(s-k)}(x)| \\ &\leq Const. \, s! \left(\sum_{k=0}^s \frac{p^{s-k}}{(s-k)!}(\eta_c|x^2 - \pi^2|)^{-k}\right) \cdot e^{\left(\frac{cx^2}{x^2 - \pi^2}\right)} \\ &\leq Const. \, \frac{s!}{(\lambda_c|x^2 - \pi^2|)^s} e^{\left(p\lambda_c|x^2 - \pi^2| + \frac{cx^2}{x^2 - \pi^2}\right)} \end{aligned}$$

which after dilation satisfies

$$|\psi_{p,\theta}^{(s)}(x)| \leq Const. \, s! \left(\frac{\theta}{\lambda_c |a(x)|}\right)^s \cdot e^{\left(\frac{p\lambda_c|a(x)|}{\theta^2} + \frac{cx^2}{a(x)}\right)}, \quad a(x) := x^2 - \pi^2\theta^2(x). \tag{3.12}$$

Following a similar manipulation we used earlier, the upper bound on right hand side of the (3.12) is maximized at $x = x_{max}$ with $x_{max}^2 - \pi^2\theta^2 \sim -c\pi^2\theta^2/s$, which leads to the $G_2$-regularity bound for $\psi_{p,\theta}$ (where as before, $\eta_c := \lambda_c \pi^2 c$)

$$\sup_{x \in \Re} |\psi_{p,\theta}^{(s)}(x)| \leq Const \cdot s! \left(\frac{s}{\eta_c \theta e}\right)^s e^{p\eta_c/s} \leq Cosnt. \frac{(s!)^2}{(\eta_c \theta)^s} e^{p\eta_c/s} \quad s = 1, 2, \ldots \tag{3.13}$$

With (3.13) we utilize (2.22) to obtain the following precise bound of the truncation error

$$|T(N, p, \theta)| \leq Const. \frac{(s!)^2}{(\eta_c \theta N)^s} e^{p\eta_c/s} \tag{3.14}$$

To minimize the upperbound (3.14), we first seek the minimizer for the order of smoothness, $s = s_p$, and then optimize the free spectral parameter $p \leq N$ for both the truncation and regularizations errors. We begin by noting that a general expression of the type encountered on the right of (3.14),

$$\frac{(s!)^2}{(\eta_c \theta N)^s} e^{p\eta_c/s} \sim \left(\frac{s^2}{\eta_c \theta e^2 N}\right)^s e^{p\eta_c/s} =: M(s, p),$$

is minimized at the $p$-dependent index $s_{min}$ such that

$$\partial_s (\log M(s, p))_{|s=s_{min}} = \log\left(\frac{s_{min}^2}{\eta_c \theta N}\right) - \frac{p\eta_c}{s_{min}^2} = 0 \tag{3.15}$$

Though we cannot find its explicit solution, (3.15) yields a rather precise bound on $s_{min}$ which turns out to be essentially independent of $p$. Indeed, for the first expression on the right of (3.15) to be positive we need $s_{min} = \sqrt{\beta \eta_c \theta N}$ with some $\beta > 1$. Plugging this expression, $s_{min} = \sqrt{\beta \eta_c \theta N}$, into (3.15), we find that for $p \leq N$ we must have, $\log \beta = \log(s_{min}^2/\eta_c \theta N) = \eta_c p/s_{min}^2 \leq 1/\beta\theta$. We therefore set $s \sim s_{min}$ of the form

$$s = \sqrt{\beta \eta_c \theta N}, \quad 1 < \beta < 1.764,$$

so that the free $\beta$ parameter satisfies the above constraint[5] $\beta \log \beta \leq 1 \leq 1/\theta$. The corresponding optimal parameter $p$ is then given by

$$p_{min} = \frac{s^2}{\eta_c} \cdot (\log \frac{s^2}{\eta_c \theta N})_{|s=s_{min}} = \kappa \cdot \theta N, \quad 0 < \kappa =: \beta \log \beta < 1 (\leq \frac{1}{\theta}). \tag{3.16}$$

We conclude with an optimal choice of $p$ of order $\mathcal{O}(\theta N)$, replacing the previous choice, (2.37), of order $\mathcal{O}(\sqrt{N})$. The resulting exponentially small truncation error bound, (3.14), now reads

$$|T(N, p, \theta)| \leq Const. \frac{(s!)^2}{(\eta_c \theta N)^s} e^{p\eta_c/s}{}_{|s=s_{min}} \sim \sqrt{\theta N} \left(\frac{\beta}{e}\right)^{2\sqrt{\beta \eta_c \theta N}}, \quad 1 < \beta \leq 1.764. \tag{3.17}$$

With this choice of $p = p_{min}$ in (3.16) we find essentially the same exponentially small bound on the regularization error in (3.10),

---

[5]Recall that $\theta = \theta(x) := d(x)/\pi < 1$.

$$|R(N, p, \theta)| \leq Const \cdot \theta N \left(\frac{1}{e}\right)^{2\sqrt{\beta \log \beta \cdot \eta_c \theta N}}. \tag{3.18}$$

Figures 3.1(f) and 3.2(f) below confirm that the contributions of the truncation and regularization parts of the error are of the same exponentially small order *up to* the vicinity of the discontinuous jumps with this choice of optimal $p \sim Nd(x)/\pi$, in contrast to previous choices of $p = \mathcal{O}(N^\gamma)$, $\gamma < 1$, consult Figures 3.1(b)-(d) and 3.2(b)-(d).

We summarize what we have shown in the following theorem.

**Theorem 3.1** *Given the Fourier projection, $S_N f(\cdot)$ of a piecewise analytic $f(\cdot)$, we consider the 2-parameter family of spectral mollifiers*

$$\psi_{p,\theta}(x) := \frac{1}{\theta}\rho_c(\frac{x}{\theta})D_p(\frac{x}{\theta}), \quad \rho_c := e^{\left(\frac{cx^2}{x^2 - \pi^2}\right)} 1_{[-\pi,\pi]}, \ c > 0,$$

*and we set*

$$\theta = \theta(x) := \frac{d(x)}{\pi}, \quad d(x) = dist(x, sing\ suppf) \tag{3.19}$$

$$p = p(x) \sim \kappa \cdot \theta(x)N, \quad 0 < \kappa = \beta \log \beta < 1. \tag{3.20}$$

*Then there exist constants, $Const_c$ and $\eta_c$, depending solely on the analytic behavior of $f(\cdot)$ in the neighborhood of $x$, such that we can recover the intermediate values of $f(x)$ with the following exponential accuracy*

$$|\psi_{p,\theta} \star S_N f(x) - f(x)| \leq Const_c \cdot \theta N \left(\frac{\beta}{e}\right)^{2\sqrt{\kappa \eta_c \theta(x) N}}, \quad 1 < \beta \leq 1.764. \tag{3.21}$$

*Remark.* Theorem 3.1 indicates an optimal choice for the spectral mollifier, $\psi_{p,\theta}$, based on an *adaptive* degree of order $p = \kappa \theta(x) N$, with an arbitrary free parameter, $0 < \kappa = \beta \log \beta < 1$. We could further optimize the error bound (3.21) over all possible choices of $\beta$, by equilibrating the leading term in the truncation and regularization error bounds so that

$$\left(\frac{\beta}{e}\right)^{2\sqrt{\beta \cdot \eta_c \theta(x) N}} \sim \left(\frac{1}{e}\right)^{2\sqrt{\beta \log \beta \cdot \eta_c \theta(x) N}},$$

with the minimal value found at $\log \beta_* = (3 - \sqrt{5})/2$, the corresponding $\kappa_* := \beta_* \log \beta_* = 0.5596$, and $2\sqrt{\kappa_*/\pi} = 0.8445$ leading to an error bound,

$$|\psi_{p,\theta} \star S_N f(x) - f(x)| \leq Const_c \cdot d(x) N \left(\frac{1}{e}\right)^{0.845\sqrt{\eta_c d(x) N}}. \tag{3.22}$$

Although the last estimate serves only as an upperbound for the error, it is still remarkable that the (close to) optimal parameterization of the adaptive mollifier is found to be essentially independent of the properties of $f(\cdot)$.

Similar result holds in the pseudospectral case. In this case, we are given the Fourier interpolant, $I_N f(x)$ and the corresponding discrete convolution is carried out in the physical space with overall error, $E(N, p, \theta; f) = \psi_{p,\theta} \star I_N f(x) - f(x)$, which consists of aliasing and regularization errors, (2.31). According to (2.32), the former is upper bounded by the truncation of $f'$, which retains the same analyticity properties as $f$ does. We conclude

**Theorem 3.2** *Given the equidistant gridvalues, $\{f(x_\nu)\}_{0 \leq \nu \leq 2N-1}$ of a piecewise analytic $f(\cdot)$, we consider the 2-parameter family of spectral mollifiers*

$$\psi_{p,\theta}(x) := \frac{1}{\theta}\rho_c(\frac{x}{\theta})D_p(\frac{x}{\theta}), \quad \rho_c := e^{\left(\frac{cx^2}{x^2 - \pi^2}\right)} 1_{[-\pi,\pi]}, \ c > 0,$$

*and we set*

$$\theta = \theta(x) := \frac{d(x)}{\pi}, \quad d(x) = dist(x, sing\ suppf) \qquad (3.23)$$

$$p = p(x) \sim \kappa \cdot \theta(x)N, \quad 0 < \kappa = \beta \log \beta < 1. \qquad (3.24)$$

*Then, there exist constants, $Const_c$ and $\eta_c$, depending solely on the analytic behavior of $f(\cdot)$ in the neighborhood of $x$, such that we can recover the intermediate values of $f(x)$ with the following exponential accuracy*

$$\left|\frac{\pi}{N}\sum_{\nu=0}^{2N-1}\psi_{p,\theta}(x-y_\nu)f(y_\nu) - f(x)\right| \leq Const_c \cdot (d(x)N)^2 \left(\frac{\beta}{e}\right)^{2\sqrt{\beta\eta_c\theta(x)N}}, \quad 1 < \beta < 1.764. \quad (3.25)$$

## 3.2 Numerical experiments

The first set of numerical experiments compares our results with those of Gottlieb-Tadmor, [GoTa85], involving the same choice of $f(\cdot) = f_1(\cdot)$

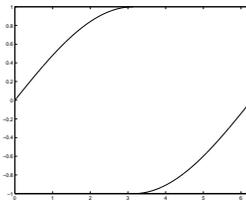

$$f_1(x) = \begin{cases} \sin(x/2) & x \in [0,\pi) \\ -\sin(x/2) & x \in [\pi,2\pi) \end{cases} \qquad (3.26)$$

A second set of results is demonstrated with a second function, $f_2(x)$ given by

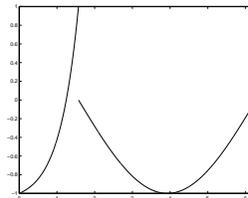

$$f_2(x) = \begin{cases} (2e^{2x} - 1 - e^\pi)/(e^\pi - 1) & x \in [0,\pi/2) \\ -\sin(2x/3 - \pi/3) & x \in [\pi/2,2\pi). \end{cases} \qquad (3.27)$$

This is a considerable challenging test problem: $f_2$ has a jump discontinuity at $x = \pi/2$ and due to the periodic extention of the Fourier series two more discontinuities at the boundaries $x = 0, 2\pi$. Moreover, a relatively large gradient is formed for $x \sim \pi/2-$, and the sharp peak on the left of $x = \pi/2$ is met by a jump discontinuity on the right.

For the computations below we utilize the same localizer $\rho = \rho_c$ as in (3.2), with $c = 10$. In the first case, $f_1$ has a simple discontinuity at $x = \pi$ so the $\theta$ parameter was chosen according to (2.36), $\theta = \theta(x) = \min(|x|, |x - \pi|)/\pi$. In the second case of $f_2(x)$ we set

$$\theta(x) = [\min(x, \pi/2 - x)_+ + \min(x - \pi/2, 2\pi - x)_+]/\pi.$$

Since the error deteriorates in the immediate vicinity of the discontinuities where $\theta(x)N \sim 1$, a window of minimum width of $\theta_{min} = \min\{\theta(x), 1/4N\}$ was imposed around $x_0 = \{0, \pi/2, \pi, 2\pi\}$. More about the treatment in the immediate vicinity of the discontinuity is found in §4.1.

The different policies for choosing the parameter $p$ are outlined below. In particular, for the near optimal choice recommended in Theorems 3.1 and 3.2 we use a mollifier of an adaptive degree $p = \kappa \theta(x) N$ with $\kappa = 1/\sqrt{e} = 0.6095 \sim \kappa_*$.

We begin with the results based on the spectral projections, $S_N f_1$ and $S_N f_2$. For comparison purposes, the exact convolution integral, $\psi_{p,\theta} \star S_N f$ was computed with composite Simpson method using $\frac{\pi}{8000}$ points, and the mollified results are recorded at the left half points, $\frac{\pi}{150}\nu$, $\nu = 0, 1, \ldots, 149$. Figure 3.1 shows the result of treating the spectral projection $\psi_{p,\theta} \star S_N f_1$ based on $N = 128$ modes, for different choices of $p$'s. Figures 3.2 show the same results for $f_2(x)$. It is evident from these figures, figure 3.1(e)-(f) and figure 3.2(e)-(f), that best results are obtained with $p = \theta(x)N/\sqrt{e}$, in agreement with our analysis for the optimal choice of exponentially accurate mollifier in (3.22). We note that other choices for $p \sim N^\gamma$, lead to large intervals where exponential accuracy is lost due to the imbalance between the truncation and regularization errors, consult cases (a)-(d) in figures 3.1-3.2. As we noted earlier in (2.38), a nonadaptive choice of $p$ independent of $\theta(x)$ leads to deterioration of the accuracy in an increasing region of size $\sim p\pi/N$ around the discontinuity, and the predicted locations of these values, given in table 3.1, could be observed in figures 3.1(a)-(d) for the function $f_1(x)$.

| $p \backslash N$ | 32 | 64 | 128 |
|---|---|---|---|
| $N^{0.8}$ | 1.6 | 1.8 | 2.0 |
| $N^{0.5}$ | 2.6 | 2.7 | 2.9 |
| $N^{0.2}$ | 2.9 | 3.0 | 3.1 |

Table 3.1: Predicted location where spectral convergence is lost at $|x - x_0| \sim p\pi/N$.

Figure 3.3 illustrates the spectral convergence as $N$ doubles from 32 to 64, then to $128^6$. The exponential convergence of the near optimal adaptive $p = \theta(x)N/\sqrt{e}$ can be seen in figure 3.3(e)-(f), where the log-slopes are constants with respect to $d(x)$ (for fixed $N$) and with respect to $N$ (for fixed $x$).

Next, the numerical experiments are repeated for the discrete case, using discrete mollification of the Fourier interpolant. Given the gridvalues of $f_1(x_\nu)$ and $f_2(x_\nu)$ at the equidistant gridpoints $x_\nu = \nu\pi/N$, we recover the pointvalues at the intermediate gridpoints $f(x_{\nu+1/2})$. A minimal window width of $\theta_{min} = \min(\theta(x), 2\pi/N)$ was imposed in the immediate vicinity of the discontinuities to maintain a minimum number of two sampling points to be used in the discrete mollification.

Compared with the previous mollified results of the spectral projections, there are two noticeable changes, both involving the non-optimal choice of $p \sim N^\gamma$ with $\gamma < 1$: (i) The location where spectral/exponential convergence is lost is noticeably closer to the discontinuities compared with the mollified spectral projections, but at the same time (ii) Much larger oscillations are observed in the regions where spectral convergence is lost. Comparing figures 3.1 vs. 3.4, and 3.2 vs. 3.5, gives a visual comparison for both changes from spectral to pseudo spectral. The deterioration for $p = N^{0.8}$ and $N^{0.5}$ are very noticeable.

---

[6]Machine truncation error is at $-16$.

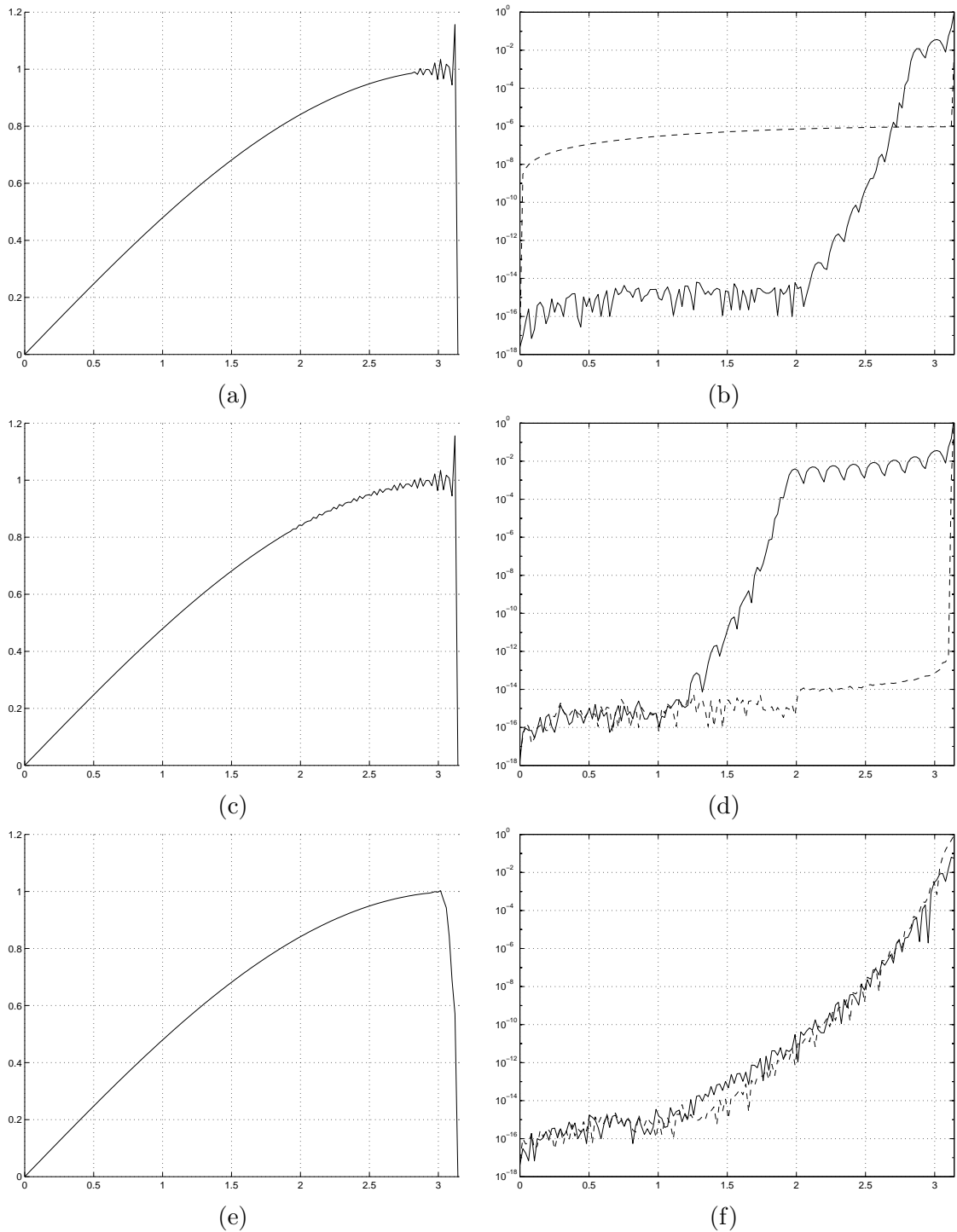

Figure 3.1: Recovery of $f_1(x)$ from its first $N = 128$ Fourier modes, on the left, and the corresponding regularization errors (dashed) and truncation errors (solid) on the right, using the spectral mollifier $\psi_{p,\theta}$ based on various choices of $p$: (a)-(b) $p = N^{0.5}$, (c)-(d) $p = N^{0.8}$, (e)-(f) $p = Nd(x)/\pi\sqrt{e}$.

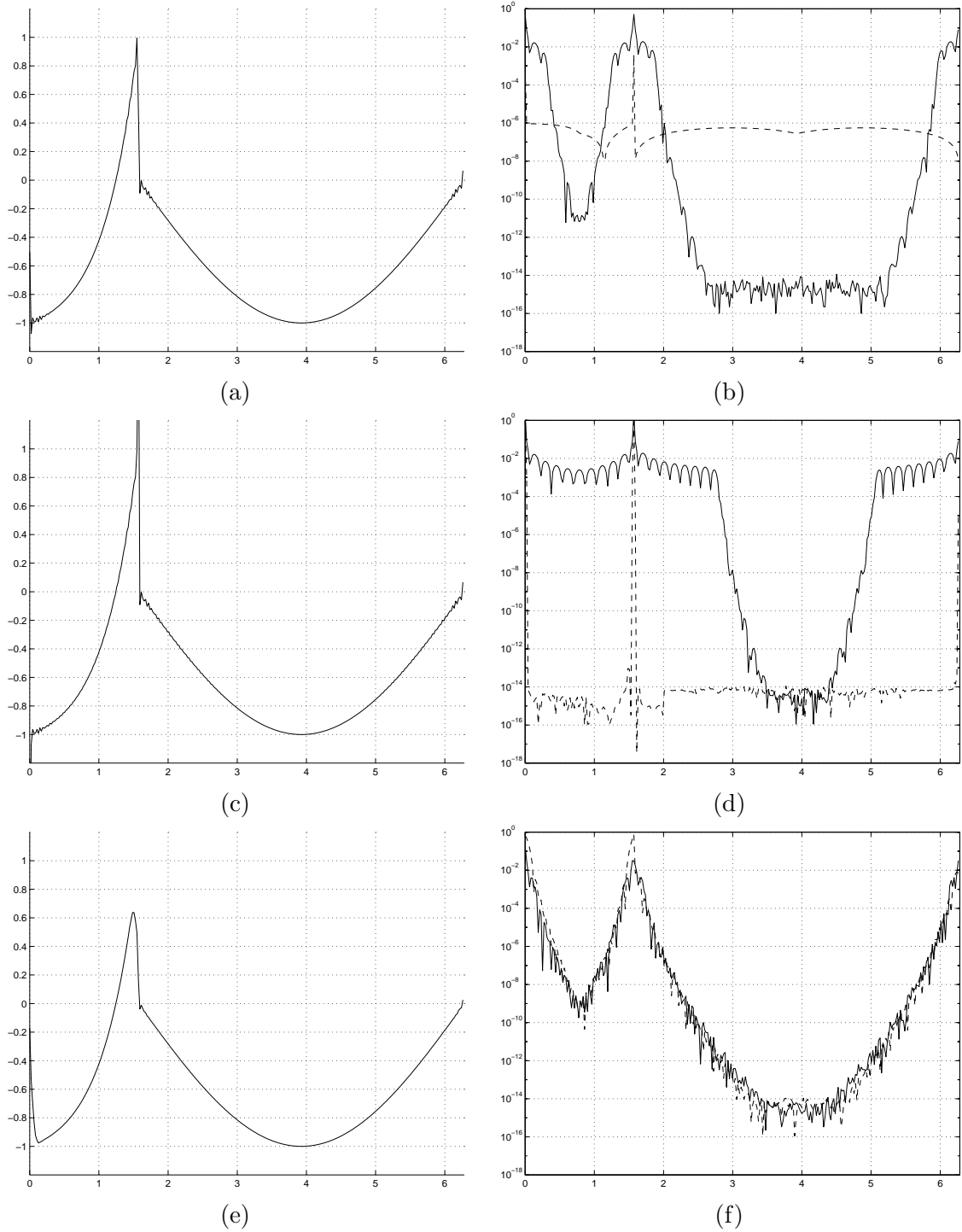

Figure 3.2: Recovery of $f_2(x)$ from its first $N = 128$ Fourier modes, on the left, and the corresponding regularization errors (dashed) and truncation errors (solid) on the right, using the spectral mollifier $\psi_{p,\theta}$ based on various choices of $p$: (a)-(b) $p = N^{0.5}$, (c)-(d) $p = N^{0.8}$, (e)-(f) $p = Nd(x)/\pi\sqrt{e}$.

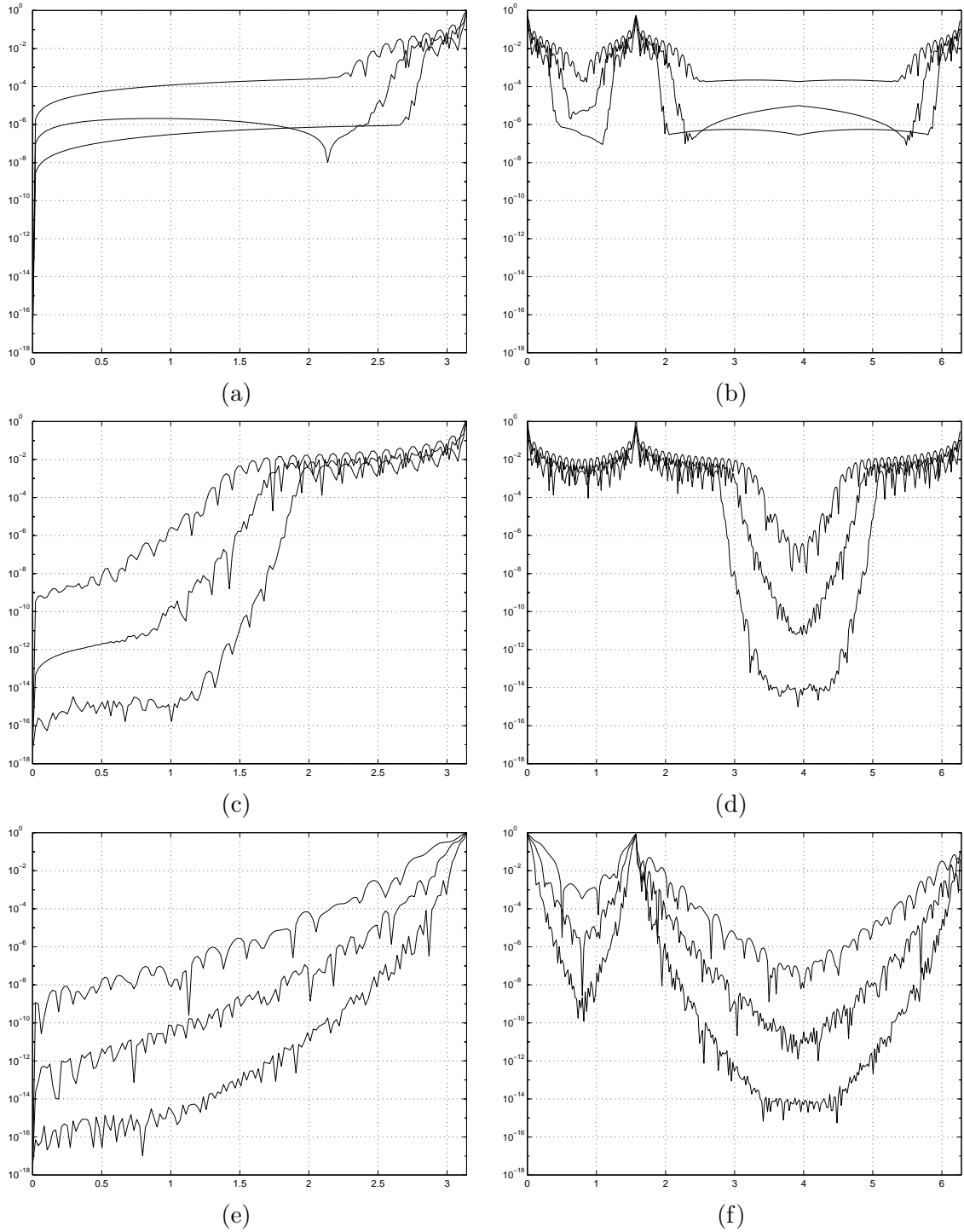

Figure 3.3: Log of the error with $N = 32, 64, 128$ modes for $f_1(x)$ on the left, and for $f_2(x)$ on the right, using various choices of $p$: (a)-(b) $p = N^{0.5}$, (c)-(d) $p = N^{0.8}$, (e)-(f) $p = Nd(x)/\pi\sqrt{e}$.

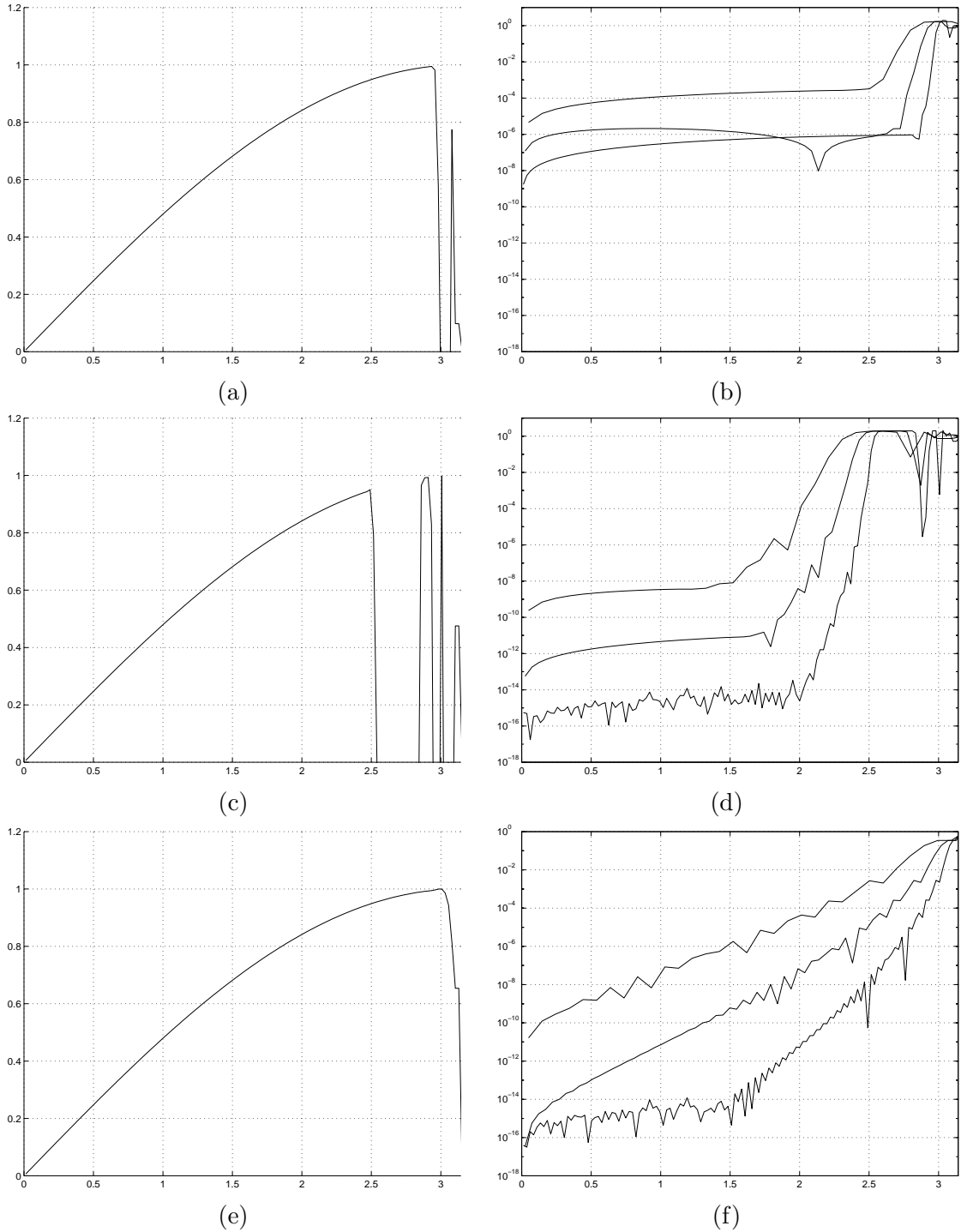

Figure 3.4: Recovery of $f_1(x)$ from its $N = 128$ equidistant gridvalues on the left, and the corresponding regularization errors (dashed) and truncation errors (solid) on the right, using the spectral mollifier $\psi_{p,\theta}$ based on various choices of $p$: (a)-(b) $p = N^{0.5}$, (c)-(d) $p = N^{0.8}$, (e)-(f) $p = Nd(x)/\pi\sqrt{e}$.

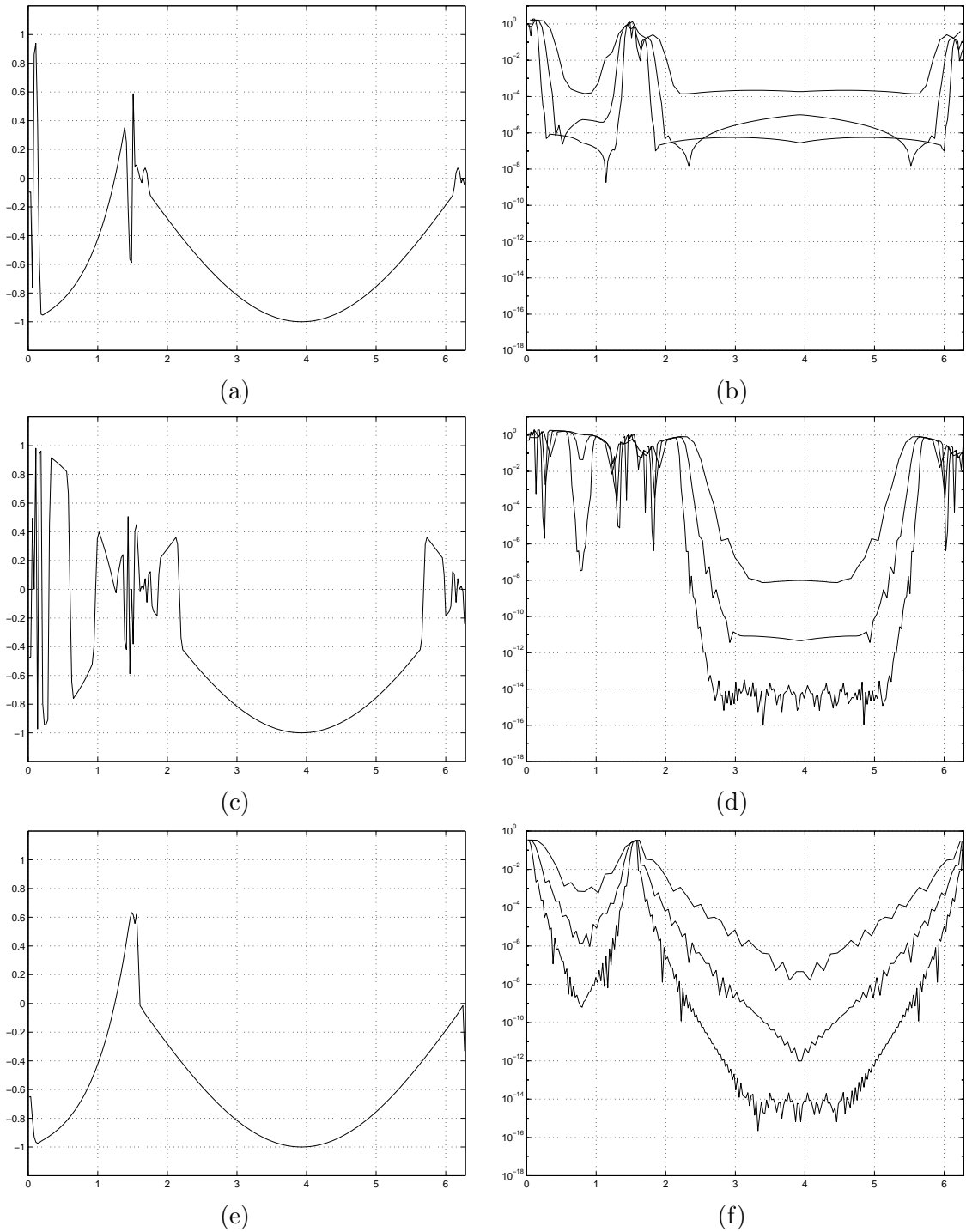

Figure 3.5: Recovery of $f_2(x)$ from its $N = 128$ equidistant gridvalues on the left, and the corresponding regularization errors (dashed) and truncation errors (solid) on the right, using the spectral mollifier $\psi_{p,\theta}$ based on various choices of $p$: (a)-(b) $p = N^{0.5}$, (c)-(d) $p = N^{0.8}$, (e)-(f) $p = Nd(x)/\pi\sqrt{e}$.

## 4   Adaptive Mollifiers – Normalization

The essence of the 2-parameter spectral mollifier discussed in §3, $\psi_{p,\theta}(x)$, is adaptivity: it is based on a Dirichlet kernel of a *variable* degree, $p \sim \theta(x)N$, which is adapted by taking into account the location of $x$ relative to its nearest singularity, $\theta(x) = d(x)/\pi$. The resulting error estimate tells us that there exist constants, $Const, \gamma$ and $\alpha > 1$ such that one can recover a piecewise analytic $f(x)$ from its spectral or pseudospectral projections, $P_N f(\cdot)$,

$$|\psi_{p,\theta} \star P_N f(x) - f(x)| \leq Const.e^{-(\gamma d(x))N)^{1/\alpha}}.$$

The error bound on the right shows that the adaptive mollifier is exponentially accurate for all $x$'s, except for what we refer to as *the immediate vicinity* of the jump discontinuities of $f$, namely, those $x$'s where $d(x) \sim 1/N$. This should be compared with previous, non-adaptive choices for choosing the degree of $\psi_{p,\theta}$: for example, with $p \sim \sqrt{N}$ we found a loss of exponential accuracy in a zone of size $\sim 1/\sqrt{N}$ around the discontinuities of $f$. Put differently, there are $\mathcal{O}(\sqrt{N})$ 'cells' which are not accurately recovered in this case. In contrast, our adaptive mollifier is exponentially accurate at all but *finitely many* cells near the jump discontinuities. According to the error estimates (3.21), (3.22) and (3.25), convergence may fail in these cells inside the immediate vicinity of $sing\ supp f$, and indeed, spurious oscillations could be noticed in figures 3.1 and 3.2. In this section we address the question of convergence *up to* the jump discontinuities.

One possible approach is to retain a uniform exponential accuracy up to the jump discontinuities. Such an approach, developed by Gottlieb, Shu, Gelb and their co-workers is surveyed in [GoSh95, GoSh98]. It is based on Gegenbauer expansions of degree $\lambda \sim N$. Exponential accuracy is retained *uniformly* throughout each interval of smoothness of the piecewise analytic $f$. The computational of the high order Gegenbauer coefficients, however, is numerically sensitive and the parameters involved need to be properly tuned in order to avoid triggering of instabilities, [Ge97, Ge00].

Here we proceed with another approach where we retain a *variable* order of accuracy near the jump discontinuities, of order $\mathcal{O}((d(x))^{r+1})$. Comparing this polynomial error bound against the interior exponential error bounds, say (3.22),

$$Const.d(x)N \cdot e^{-0.845\sqrt{\eta_c d(x)N}} \geq (d(x))^{r+1},$$

we find that there are only *finitely many* cells in which the error – dictated by the smaller of the two, is dominant by polynomial accuracy

$$d(x) \leq Const.\frac{r^2(\log d(x))^2}{\eta_c N} \sim \frac{r^2}{N}. \tag{4.1}$$

In this approach, the variable order of accuracy suggested by (4.1), $r \sim \sqrt{d(x)N}$, is increasing together with the increasing distance away from the jumps, or more precisely – together with the number of cells away from the discontinuities, which is consistent with the *adaptive* nature of our exponentially accurate mollifier away from the immediate vicinity of these jumps. The current approach of variable order of accuracy which adapted to the distance from the jump discontinuities, is reminiscent of the Essentially Non Oscillatory (ENO) piecewise polynomial reconstruction employed in the context of nonlinear conservation laws [HEOC85, Sh97].

How to enforce that our adaptive mollifiers are polynomial accurate in the immediate vicinity of jump discontinuities? as we argued earlier in (2.11), the adaptive mollifier $\psi_{p,\theta}$ admits spectrally small moments of order $p^{-s} \sim (d(x)N)^{-s}$. More precisely, using (3.7) we find for $\rho \in G_\alpha$,

$$\int_{-\pi\theta}^{\pi\theta} y^s \psi_{p,\theta}(y) dy \;=\; \int_{-\pi}^{\pi} (y\theta)^s \rho(y) D_p(y) dy = D_p \star \Big((y\theta)^s \rho\Big)(y)_{|y=0} =$$

$$= \delta_{s0} + Const.\theta^s \cdot e^{-\alpha(\eta p)^{1/\alpha}}, \qquad p \sim d(x)N.$$

Consequently, $\psi_{p,\theta}$ possesses exponentially small moments at all $x$'s except for the immediate vicinity of the jumps where $p \sim d(x)N \sim 1$, the same $\mathcal{O}(1/N)$ neighborhoods where the previous exponential error bounds fail. This is illustrated in the numerical experiments exhibited in §3.2 which show the blurring in symmetric intervals with width $\sim 1/N$ around each discontinuity. To remove this blurring, we will impose a novel *normalization* so that finitely many moments of (the projection of) $\psi_{p,\theta}$ *precisely* vanish. As we shall see below, this will regain a polynomial convergence rate of the corresponding finite order $r$. We have seen that the general adaptivity (4.1) requires $r \sim \sqrt{d(x)N}$; in practice, enforcing a fixed number of vanishing moments, $r \sim 2,3$ will suffice.

## 4.1 Spectral normalization - adaptive mollifiers in the vicinity of jumps

Rather than $\psi_{p,\theta}$ possessing a fixed number of vanishing moments as in standard mollification (2.7), we require that its spectral projection, $S_N \psi_{p,\theta}$, posses a unit mass and, say $r$ vanishing moments,

$$\int_{-\pi}^{\pi} y^s (S_N \psi_{p,\theta})(y) dy = \delta_{s0} \qquad s = 0, 1, \ldots, r. \tag{4.2}$$

It then follows that adaptive mollification of the Fourier projection, $\psi_{p,\theta} \star S_N f$, recovers the point-values of $f$ with the desired polynomial order $\mathcal{O}(d(x))^r$. Indeed, noting that for each $x$, the function $f(x-y)$ remains smooth in the neighborhood $|y| \leq \pi\theta = d(x)$ we find, utilizing the symmetry of the spectral projection, $\int (S_N f)g = \int f(S_N g)$,

$$\psi_{p,\theta} \star S_N f(x) - f(x) = \int_{-\pi\theta}^{\pi\theta} [f(x-y) - f(x)](S_N \psi_{p,\theta})(y) dy =$$
$$= \sum_{s=1}^{r} \frac{(-1)^s}{s!} f^{(s)}(x) \int_{-\pi}^{\pi} y^s (S_N \psi_{p,\theta})(y) dy + \frac{(-1)^{r+1}}{(r+1)!} f^{(r+1)}(\cdot) \int_{-\pi}^{\pi} y^{r+1} (S_N \psi_{p,\theta})(y) dy$$
$$\sim \int_{-\pi}^{\pi} (S_N y^{r+1}) \psi_{p,\theta}(y) dy \leq Const. \Big(d(x) + \frac{1}{N}\Big)^{r+1}.$$

The last step follows from an upperbound for the spectral projection of monomials outlined at the end of this subsection.

To enforce the vanishing moments condition (4.2) on the adaptive mollifier, $\psi_{p,\theta}(\cdot) = (\rho(\cdot) D_p(\cdot/\theta))/\theta$, we take advantage of the freedom we have in choosing the localizer $\rho(\cdot)$. We begin by normalizing

$$\widetilde{\psi}_{p,\theta}(y) = \frac{\psi_{p,\theta}(y)}{\int_{-\pi}^{\pi} \psi_{p,\theta}(z) dz}$$

so that $\widetilde{\psi}_{p,\theta}$ has a unit mass, and hence (4.2) holds for $r = 0$, for $\int S_N(\psi_{p,\theta})(y) dy = \int \psi_{p,\theta}(y) dy = 1$. We note that the resulting mollifier takes the same form as before, namely

$$\widetilde{\psi}_{p,\theta}(y) := \frac{1}{\theta}(\tilde{\rho}_c D_p)(\frac{y}{\theta}), \tag{4.3}$$

where the only difference is associated with the modified localizer,

$$\tilde{\rho}_c(y) = q_0 \cdot \rho_c(y), \qquad q_0 = \frac{1}{\int_{-\pi}^{\pi} \psi_{p,\theta}(z) dz}. \tag{4.4}$$

Observe that in fact, $1/q_0 = \int \psi_{p,\theta}(z)dz \equiv \int \psi_p(z)dz = (D_p \star \rho_c)(0)$, and that with our choice of $p = \kappa \cdot \theta(x)N$, we have in view of (3.7),

$$\tilde{\rho}_c(0) = q_0 = \frac{1}{(D_p \star \rho_c)(0)} = 1 + \mathcal{O}(\varepsilon), \qquad \varepsilon \sim d(x)N \cdot e^{-2\sqrt{\eta_c p}}, \quad p = \kappa \cdot \theta(x)N$$

which is admissible within the same exponentially small error bound we had before, consult (3.11). In other words, we are able to modify the localizer $\rho_c(\cdot) \to \tilde{\rho}_c(\cdot)$ to satisfy the first-order normalization, (4.2) with $r = 0$, while the corresponding mollifier, $(\rho_x D_p)_\theta \to (\tilde{\rho}_c D_p)_\theta$, retains the same overall exponential accuracy. Moreover, using even $\rho$'s implies that $\psi(\cdot)$ is an even function and hence its odd moments vanish. Consequently, (4.2) holds with $r = 1$, and we end up with the following quadratic error bound in the vicinity of $x$ (compared with (3.22))

$$|\widetilde{\psi}_{p,\theta}(x) \star S_N f(x) - f(x)| \leq Const. \left(d(x) + \frac{1}{N}\right)^2 \cdot e^{-0.845\sqrt{\eta_c d(x)N}}.$$

In a similar manner, we can enforce higher vanishing moments by proper *normalization* of the localizer $\rho(\cdot)$. There is clearly more than one way to proceed – here is one possibility. In order to satisfy (4.2) with $r = 2$ we use a pre-factor of the form $\tilde{\rho}_c(x) \sim (1 + q_2 x^2)\rho_c(x)$. Imposing a unit mass and vanishing second moment we may take

$$\widetilde{\psi}_{p,\theta}(y) = \frac{1}{\theta}(\tilde{\rho}_c D_p)\left(\frac{y}{\theta}\right), \qquad \tilde{\rho}_c(y) \sim (1 + q_2 y^2)\rho_c(y),$$

with the normalized localizer, $\tilde{\rho}_c(y)$, given by

$$\tilde{\rho}_c(y) = \frac{1 + q_2 y^2}{\int_{-\pi}^{\pi}(1 + q_2(\frac{z}{\theta})^2)\psi_{p,\theta}(z)dz}\rho_c(y), \qquad q_2 = \frac{-\int_{-\pi}^{\pi}(S_N z^2)\psi_{p,\theta}(z)dz}{\int_{-\pi}^{\pi}(S_N z^2)(\frac{z}{\theta})^2\psi_{p,\theta}(z)dz}. \tag{4.5}$$

As before, the resulting mollifier $\widetilde{\psi}_{p,\theta}$ is admissible in the sense of satisfying the normalization (3.11) within the exponentially small error term. Indeed, since $\int y^2 \psi_{p,\theta}(y)dy = (D_p \star (y^2 \rho_c(y)))(0) = \mathcal{O}(\varepsilon)$ we find

$$\tilde{\rho}_c(0) = \frac{1}{\int_{-\pi}^{\pi}(1 + q_2(\frac{y}{\theta})^2)\psi_{p,\theta}(y)dy} =$$
$$= \frac{1}{1 + q_2 \cdot \varepsilon/\theta^2} = 1 + Const \cdot (d(x)N)^3 \cdot e^{-2\sqrt{\eta_c p}}, \qquad p = \kappa \cdot \theta(x)N.$$

A straightforward computation shows that the unit mass $\widetilde{\psi}_{p,\theta}$ has a second vanishing moment

$$\int_{-\pi}^{\pi} y^2 (S_N \widetilde{\psi}_{p,\theta})(y)dy = \int_{-\pi}^{\pi}(S_N y^2)\left(a_0 + a_2\left(\frac{y}{\theta}\right)^2\right)\psi_{p,\theta}(y)dy =$$
$$= \int_{-\pi}^{\pi}(S_N y^2)\psi_{p,\theta}(y)dy + q_2 \int_{-\pi}^{\pi}(S_N y^2)\left(\frac{y}{\theta}\right)^2 \psi_{p,\theta}(y)dy = 0. \tag{4.6}$$

Since $\tilde{\rho}_c(\cdot)$ is even, so is the normalized mollifier $\widetilde{\psi}(\cdot)$, and hence its third moment vanishes yielding a 4th order convergence rate in the immediate vicinity of the jump discontinuities,

$$|(\widetilde{\psi}_{p,\theta}(x) \star S_N f)(x) - f(x)| \leq Const. \left(d(x) + \frac{1}{N}\right)^4 \cdot e^{-0.845\sqrt{\eta_c d(x)N}}.$$

We close this section with the promised

**Lemma 4.1** ([Tao]). *The following pointwise estimate holds*

$$|S_N(y^r)| \lesssim (|y| + \frac{1}{N})^r.$$

To prove this, we use a dyadic decomposition (similar to the Littlewood-Paley construction) to split

$$y^r = \sum_{k \leq 0} 2^{kr} \psi(y/2^k)$$

where $\psi$ is a bump function adapted to the set $\{\pi/4 < |y| < \pi\}$.
For $2^k \lesssim 1/N$, the usual upperbounds of the Dirichlet kernel tell us that

$$|S_N(\psi(\cdot/2^k))(y)| \lesssim 2^k N/(1 + N|y|).$$

Now suppose $2^k \gtrsim 1/N$. In this case we can use the rapid decay of the Fourier transform of $\psi(\cdot/2^k)$ for frequencies $\gg N$ to obtain the estimate

$$\|(1 - S_N)(\psi(\cdot/2^k))\|_\infty \lesssim (2^k N)^{-100}.$$

In particular, since $supp\, \psi \sim 1$, we have $|S_N(\psi(\cdot/2^k))(x)| \lesssim 1$ when $|y| \sim 2^k$, and $|S_N(\psi(\cdot/2^k))(y)| \lesssim (2^k N)^{-100}$ otherwise. The desired bound follows by adding together all these estimates over $k$.

## 4.2   Pseudospectral normalization – adaptive mollifiers in the vicinity of jumps

We now turn to the pseudospectral case which will only require evaluations of discrete sums and consequently, can be implemented with little increase in computation time.
Let $f * g(x) := \sum_\nu f(x - y_\nu) g(y_\nu) h$ denote the (non-commutative) discrete convolution based on $2N$ equidistant gridpoints, $y_\nu = \nu h, h = \pi/N$. Noting that for each $x$, the function $f(y)$ remains smooth in the neighborhood $|x - y| \leq \pi\theta = d(x)$, we find

$$|\psi_{p,\theta} * I_N f(x) - f(x)| = \Big|\sum_\nu \psi_{p,\theta}(x - y_\nu)[f(y_\nu) - f(x)]h\Big| =$$

$$= \Big|\sum_{s=1}^r \frac{(-1)^s}{s!} f^{(s)}(x) \sum_\nu (x - y_\nu)^s \psi_{p,\theta}(x - y_\nu)h + \frac{(-1)^{r+1}}{(r+1)!} f^{(r+1)}(\cdot) \sum_\nu (x - y_\nu)^{r+1} \psi_{p,\theta}(x - y_\nu)h\Big|$$

$$\leq Const.(d(x))^{r+1}, \qquad Const \sim \frac{\|f\|_{C_{loc}^{r+1}}}{(r+1)!},$$

provided $\psi_{p,\theta}$ has its first $r$ *discrete* moments vanish,

$$\sum_{\nu=0}^{2N-1} (x - y_\nu)^s \psi_{p,\theta}(x - y_\nu)h = \delta_{s0}, \qquad s = 0, 1, 2, \ldots, r. \tag{4.7}$$

Observe that unlike the continuous case associated with spectral projections, the discrete constraint (4.7) is not translation invariant and hence it requires $x$-dependent normalizations. The additional computational effort is minimal, however, due to the discrete summations which are localized in the immediate vicinity of $x$. Indeed, as a first step we note the validity of (4.7) for $x$'s which are away from the immediate vicinity of the jumps of $f$. To this end we apply the main exponential error estimate (3.25) for $f(\cdot) = (x - \cdot)^s$ (for arbitrary fixed $x$), to obtain

$$\sum_{\nu=0}^{2N-1}(x-y_\nu)^s \psi_{p,\theta}(x-y_\nu)h = (x-y)^s_{|y=x} + \mathcal{O}(\varepsilon) =$$

$$= \delta_{s0} + \mathcal{O}(\varepsilon), \qquad \varepsilon \sim (d(x)N)^2 \cdot e^{-\sqrt{Const.d(x)\cdot N}}. \qquad (4.8)$$

Thus, (4.7) holds modulo exponentially small error for those $x$'s which are away from the jumps of $f$, where $d(x) \gg 1/N$. The issue now is to enforce discrete vanishing moments on the adaptive mollifier $\psi_p(x) = \rho(x)D_p(x)$ in the vicinity of these jumps, and to this end we take advantage of the freedom we have in choosing the localizer $\rho(\cdot)$. We begin by normalizing

$$\widetilde{\psi}_{p,\theta}(y) = \frac{\psi_{p,\theta}(y)}{\sum_{\nu=0}^{2N-1}\psi_{p,\theta}(x-y_\nu)h},$$

so that $\widetilde{\psi}_{p,\theta}(x-\cdot)$ has a (discrete) unit mass, i.e., (4.7) holds with $r=0$. We note that the resulting mollifier takes the same form as before, namely

$$\widetilde{\psi}_{p,\theta}(y) := \frac{1}{\theta}(\tilde{\rho}_c D_p)(\frac{y}{\theta}), \qquad (4.9)$$

and that the only difference is associated with the modified localizer,

$$\tilde{\rho}_c(y) = q_0 \cdot \rho_c(y), \qquad q_0 = \frac{1}{\sum_{\nu=0}^{2N-1}\psi_{p,\theta}(x-y_\nu)h}. \qquad (4.10)$$

By (4.8), the $x$-dependent normalization factor, $q_0 = q_0(x)$ is in fact an approximate identity,

$$1/q_0 = \sum_{\nu=0}^{2N-1}\psi_{p,\theta}(x-y_\nu)h = 1 + \mathcal{O}(\varepsilon), \qquad \varepsilon \sim (d(x)N)^2 \cdot e^{-\sqrt{Const.d(x)\cdot N}},$$

which shows that the normalized localizer is admissible, $|\tilde{\rho}(0)-1| = |q_0-1| \leq \mathcal{O}(\varepsilon)$, within the same exponentially small error bound we had before – consult (3.11) with our choice of $p \sim d(x) \cdot N$. In other words, we are able to modify the localizer $\rho_c(\cdot) \to \tilde{\rho}_c(\cdot)$ to satisfy the first-order normalization, (4.7) with $r=0$ required near jump discontinuities, while the corresponding mollifier, $(\rho_x D_p)_\theta \to (\tilde{\rho}_c D_p)_\theta$, retains the same overall exponential accuracy required outside the immediate vicinity of these jumps .

Next, we turn to enforce that first discrete moment vanishes, $\sum_\nu (x-y_\nu)\widetilde{\psi}_{p,\theta}(x-y_\nu)h = 0$, and to this end we seek a modified mollifier of the form

$$\widetilde{\psi}_{p,\theta}(y) = \frac{q(y/\theta)}{\sum_\nu q(\frac{x-y_\nu}{\theta})\psi_{p,\theta}(x-y_\nu)h}\psi_{p,\theta}(y), \qquad q(y) := 1 + q_1 y,$$

with $q_1$ is chosen so that the second constraint, (4.7) with $r=1$, is satisfied[7]

$$q_1 = -\frac{\sum_\nu (x-y_\nu)\psi_{p,\theta}(x-y_\nu)h}{\sum_\nu \frac{(x-y_\nu)^2}{\theta}\psi_{p,\theta}(x-y_\nu)h}. \qquad (4.11)$$

---

[7]We note in passing that $\tilde{\rho}_c(\cdot)$ being even implies that $\widetilde{\psi}_{p,\theta}(\cdot)$ is an even function and hence its odd moments vanish. It follows that the first discrete moment, $\sum_\nu (x-y_\nu)\psi_{p,\theta}(x-y_\nu)h$ vanishes at the gridpoints $x=y_\mu$, and therefore $q_1 = 0$ there. But otherwise, unlike the similar situation with the spectral normalization, $q_1 \neq 0$. The discrete summation in $q_1$, however, involves only finitely many neighboring values in the $\theta$-vicinity of $x$.

Consequently, (4.7) holds with $r = 1$, and we end up with a quadratic error bound corresponding to (3.22)

$$|\widetilde{\psi}_{p,\theta} * I_N f(x) - f(x)| \leq Const.(d(x))^2 \cdot e^{-\sqrt{Const.d(x)N}}.$$

Moreover, (4.8) implies that $q_1 = \mathcal{O}(\varepsilon)$ and hence the new normalized localizer is admissible, $\tilde{\rho}_c(0) = 1 + \mathcal{O}(\varepsilon)$. In a similar manner we can treat higher moments, using *normalized* localizers, $\tilde{\rho}_c(y) \sim q(y)\rho_c(y)$ of the form

$$\widetilde{\psi}_{p,\theta}(y) = \frac{1}{\theta}(\tilde{\rho}_c D_p)(\frac{y}{\theta}), \qquad \tilde{\rho}_c(y) := \frac{1 + q_1 y + \ldots q_r y^r}{\sum_\nu q(\frac{x-y_\nu}{\theta})\psi_{p,\theta}(x - y_\nu)h}\rho_c(y). \qquad (4.12)$$

The $r$ free coefficients of $q(y) = 1 + q_1 y + \ldots q_r y^r$ are chosen so as to enforce (4.7) with the first $r$ discrete moments of $\widetilde{\psi}$ vanish. This leads to a simple $r \times r$ Vandermonde system ( – outlined in at the end of this section) involving the $r$ gridvalues, $\{f(y_\nu)\}$, in the vicinity of $x$, $|y_\nu - x| \leq \theta(x)\pi$. With our choice of a symmetric support of size $\theta(x) = d(x)/\pi$, there are precisely $r = 2\theta\pi/h = 2Nd(x)/\pi$ such gridpoints in the immediate vicinity of $x$, which enable us to recover the intermediate gridvalues, $f(x)$ with an *adaptive* order $(d(x))^{r+1}$, $r \sim Nd(x)$. As before, this normalization does not affect the exponential accuracy away from the jump discontinuities, noting that $\tilde{\rho}(0)_c = 1/q_0 = 1 + \mathcal{O}(\varepsilon)$ in agreement with (3.11). We summarize by stating

**Theorem 4.1** *Given the equidistant gridvalues, $\{f(x_\nu)\}_{0 \leq \nu \leq 2N-1}$ of a piecewise analytic $f(\cdot)$, we want to recover the intermediate values $f(x)$. To this end, we use the 2-parameter family of pseudospectral mollifiers*

$$\widetilde{\psi}_{p,\theta}(y) := \frac{1}{\theta}\tilde{\rho}_c(\frac{y}{\theta})D_p(\frac{y}{\theta}), \quad p = 0.5596 \cdot \theta N, \; c > 0,$$

*where $\theta = \theta(x) := d(x)/\pi$ is the (scaled) distance between $x$ and its nearest jump discontinuity. We set $\tilde{\rho}_c(y) := q(y)e^{\left(\frac{cy^2}{y^2-\pi^2}\right)}1_{[-\pi,\pi]}$ as the normalizing factor, with*

$$q(y) = \frac{1 + q_1 y + \ldots q_r y^r}{\sum_\nu q(\frac{x-y_\nu}{\theta})\psi_{p,\theta}(x - y_\nu)h}$$

*so that the first $r$ discrete moments of $\widetilde{\psi}_{p,\theta}(y)$ vanish, i.e., (4.7) holds with $r \sim Nd(x)$.*
*Then, there exist constants, $Const_c$ and $\eta_c$, depending solely on the analytic behavior of $f(\cdot)$ in the neighborhood of $x$, such that we can recover the intermediate values of $f(x)$ with the following exponential accuracy*

$$\left|\frac{\pi}{N}\sum_{\nu=0}^{2N-1}\psi_{p,\theta}(x - y_\nu)f(y_\nu) - f(x)\right| \leq Const_c \cdot (d(x))^{r+1}\left(\frac{1}{e}\right)^{0.845\sqrt{\eta_c d(x)N}}, \quad r \sim Nd(x). \quad (4.13)$$

The error bound (4.13) confirms our statement in the introduction of §4, namely, the *adaptivity* of the spectral mollifier in the sense of recovering the gridvalues in the vicinity of the jumps with an increasing order, $Nd(x)$, proportional to their distance from *sing supp f*. We have seen that the general adaptivity (4.1) requires $r \sim \sqrt{d(x)N}$, so that in practice, enforcing a fixed number of vanishing moments, $r \sim 2, 3$ will suffice as a transition to the exponentially error decay in the interior region of smoothness. We highlight the fact that the modified mollifier $\widetilde{\psi}_{p,\theta}$ normalized by having finitely many ($\sim 2, 3$) vanishing moments can be constructed with little increase in computation time and, as we will see in §4.3 below, it yields greatly improved results near the discontinuities.

We close this section with a brief outline on the construction of the $r$-order accurate normalization factor $q(\cdot)$. To recover $f(x)$, we seek a $r$-degree polynomial $q(y) := 1 + q_1 y + \ldots + q_r y^r$ so that (4.7) holds. We emphasize that the $q_r$'s depend on the specific point $x$ in the following manner. Setting $z_\nu := x - y_\nu$, then satisfying (4.7) for the *higher* moments of $\widetilde{\psi}_{p,\theta}$ requires

$$\sum_\nu z_\nu^s \widetilde{\psi}_{p,\theta}(z_\nu) h = 0, \qquad s = 1, 2, \ldots, r,$$

and with $\widetilde{\psi}_{p,\theta}(\cdot) \sim q(\cdot/\theta)\psi_{p,\theta}(\cdot)$ we end up with

$$\sum_\nu z_\nu^s [q(\frac{z_\nu}{\theta}) - 1]\psi_{p,\theta}(z_\nu) h = -\sum_\nu z_\nu^s \psi_{p,\theta}(z_\nu), \qquad s = 1, 2, \ldots, r.$$

Expressed in terms of the discrete moments of $\psi$,

$$a_\alpha(z_\nu) := \sum_\nu (\frac{z_\nu}{\theta})^{1+\alpha} \psi_{p,\theta}(z_\nu), \qquad \alpha = 1, 2 \ldots, 2r$$

this amounts to the $r \times r$ Vandermonde-like system for $\{q_1, \ldots, q_r\}$,

$$\begin{bmatrix} a_1(z_\nu) & a_2(z_\nu) & \cdots & a_r(z_\nu) \\ \cdot & \cdot & \cdots & \cdot \\ \cdot & \cdot & \cdots & \cdot \\ \cdot & \cdot & \cdots & \cdot \\ a_{r+1}(z_\nu) & a_{r+2}(z_\nu) & \cdots & a_{2r}(z_\nu) \end{bmatrix} \begin{bmatrix} q_1 \\ \cdot \\ \cdot \\ \cdot \\ q_r \end{bmatrix} = -\begin{bmatrix} \sum_\nu z_\nu \psi_{p,\theta}(z_\nu) \\ \cdot \\ \cdot \\ \cdot \\ \sum_\nu z_\nu^r \psi_{p,\theta}(z_\nu) \end{bmatrix} \qquad (4.14)$$

Finally we scale $q(\cdot)$ so that (4.7) holds with $s = 0$, which led us to the normalized localizer in (4.12).

### 4.3    Numerical Experiments

Figure 3.1 (d) the blurring oscillations near the edges when using the non-normalized adaptive mollifier. To reduce this blurring we will use the normalized $\widetilde{\psi}_{p,\theta}$ for $x$'s in the vicinity of the jumps where $d(x) \leq 6\pi/N$. The convolution is computed at the same locations as in section (3.2), and a minimum window width of $\theta(x) = \min(d(x)/\pi, 2\pi/N)$ was imposed. the Trapezoidal rule (with spacing of $\pi/8000$) was used for the numerical integration of $(S_N y^2)\psi_{p,\theta}(y)$ and $(S_N y^2)y^2\psi_{p,\theta}(y)$, required for the computation of $q_0$ and $q_2$ in 4.1. Figure 4.1(a)-(d) shows the clear improvement near the edges once we utilize the normalized $\widetilde{\psi}_{p,\theta}$, while retaining the exponential convergence away from these edges is illustrated in figure 4.1(e)-(f).

We conclude with the pseudospectral case. The $\mathcal{O}(1)$ error remains in figure 3.4 (d) for the non-normalized mollifier. The normalization of the discrete mollifier in section 4.2 shows that by using $\widetilde{\psi}_{p,\theta}$ given in (4.12), with a 4th degree normalization factor $q(\cdot)$, results in a minimum convergence rate of $d(x)^4$ in the vicinity of the jumps, and with exponentially increasing order as we move away from the jumps. This modification of $\widetilde{\psi}_{p,\theta}$ leads to a considerable improvement in the resolution near the discontinuity, which could be seen in figure 4.2. Here, normalization was implemented using $\widetilde{\psi}_{p,\theta}$ in the vicinity of the jumps, for $d(x) \leq 4\pi/N$, and the adaptive mollifier $\psi_{p,\theta}$ was used for $x$'s 'away' from the jumps $d(x) \geq 4\pi/N$. A minimum window of width $\theta(x) = min(d(x)/\pi, 2\pi/N)$ was imposed.

*

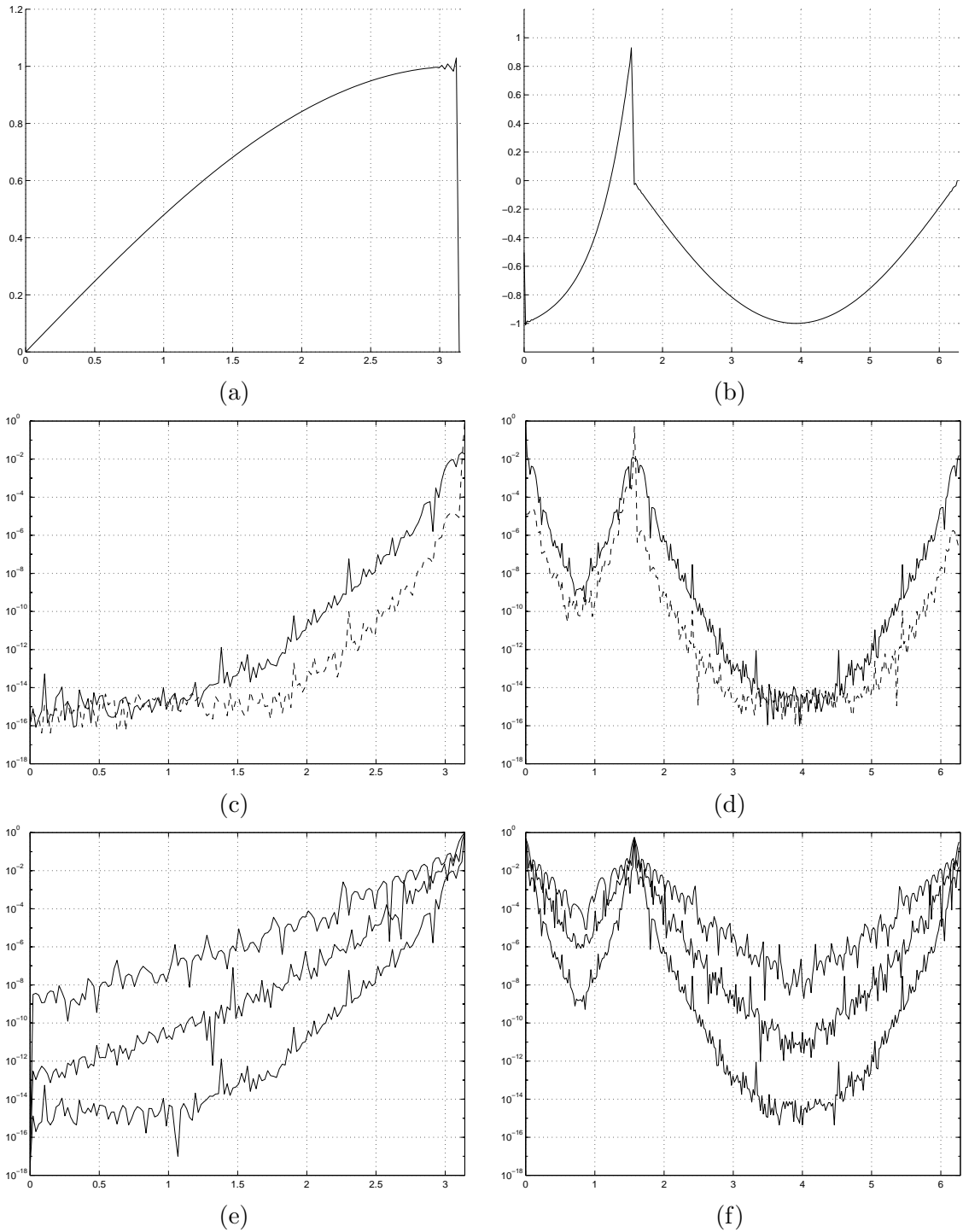

Figure 4.1: Recovery of $f_1(x)$ (on the left) and $f_2(x)$ (on the right) from their $N = 128$-modes spectral projections, using the 4th order normalized mollifier (4.3),(4.5) of degree $p = d(x)N/\pi\sqrt{e}$. Regularization errors (dashed) and truncation errors (solid) are shown on (c)-(d), and Log errors based on $N = 32, \ 64,$ and $128$ modes are shown in (e)-(f).

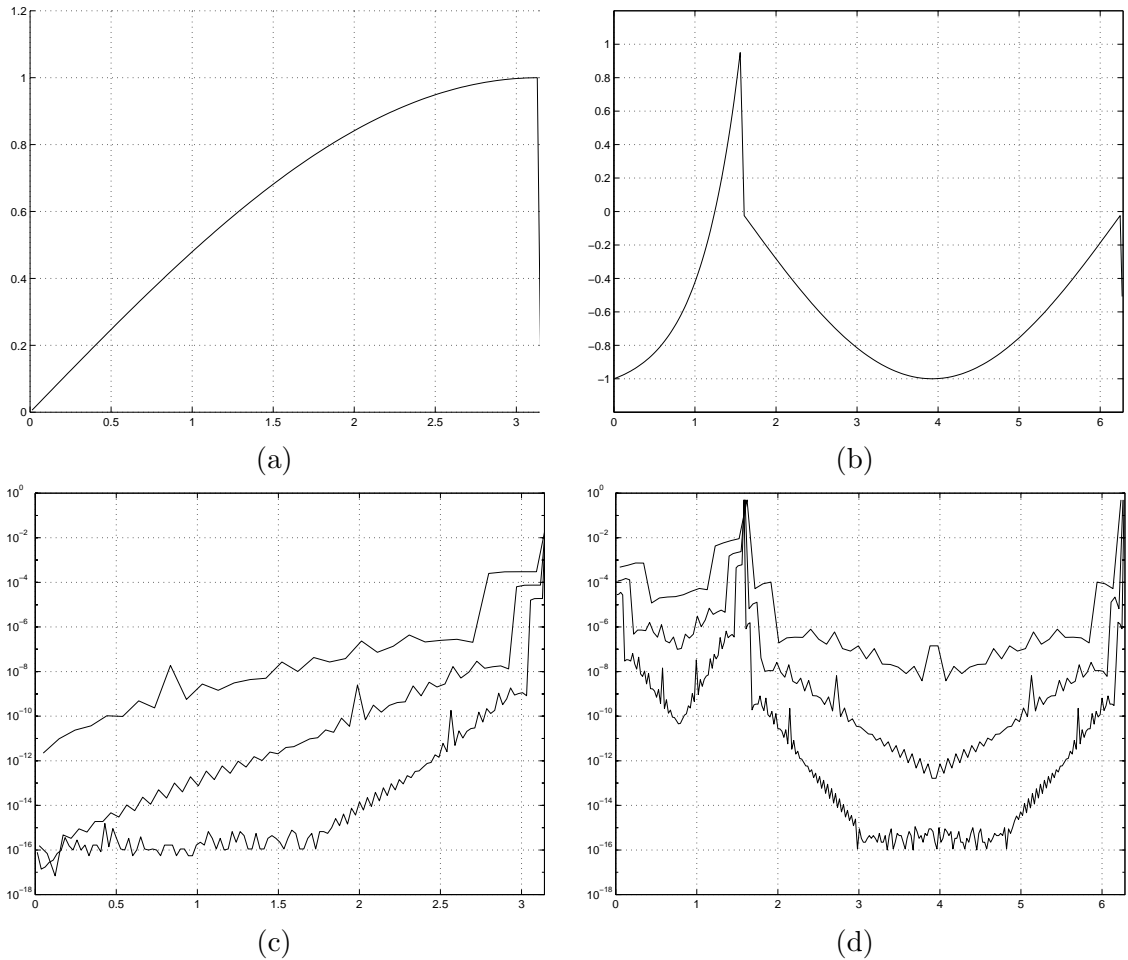

Figure 4.2: Recovery of $f_1(x)$ (a) and $f_2(x)$ (b) from their $N = 128$-modes spectral projections, using the normalized mollifier. Log error for recovery of $f_1(x)$ (c) and $f_2(x)$ (d) from their spectral projections based on $N = 32$, $64$, and $128$ modes. Here we use the normalized mollifier, $\psi_{p,\theta}$ of degree $p = d(x)N/\pi\sqrt{e}$.

## 5  Summary

In their original work [GoTa85], Gottlieb & Tadmor showed how to regain *formal* spectral convergence in recovering piecewise smooth functions using the 2-parameter family of mollifiers $\psi_{p,\theta}$. Our analysis shows that with a proper choice of parameters – in particular, an *adaptive* choice for the degree $p \sim d(x)N$, hide the overall strength in the method. By incorporating the distance to the discontinuities, $\theta = d(x)/\pi$ along with the optimal value of $p$, we end up with exponentially accurate recovery procedure up to the immediate vicinity of the jump discontinuities. Moreover, with a proper local *normalization* of the spectral mollifier, one can further reduce the error in the vicinity of these jumps. For the pseudospectral case, the normalization adds little to the overall computation time. Overall, this yields a high resolution yet very robust recovery procedure which enables one to effectively manipulate pointwise values of piecewise smooth data.

# References

[BL93] R. K. Beatson and W. A. Light, *Quasi-interpolation by thin-plate splines on the square*, Constr. Approx., 9 (1993) 343-372.

[Ch] E. W. Cheney, *Approximation Theory*, Chelsea, 1982.

[Ge97] A. Gelb, *The resolution of the Gibbs phenomenon for spherical harmonics*, Math. Comp., 66 (1997), 699-717.

[Ge00] A. Gelb, *A Hybrid Approach to Spectral Reconstruction of Piecewise Smooth Functions*, Journal of Scientific Computing, October 2000.

[GeTa99] A. Gelb and E. Tadmor, *Detection of Edges in Spectral Data*, Applied Computational Harmonic Analysis 7, (1999) 101-135.

[GeTa00a] A. Gelb and E. Tadmor, *Enhanced spectral viscosity approximations for conservation laws*, Applied Numerical Mathematics 33 (2000), 3-21.

[GeTa00b] A. Gelb and E. Tadmor, *Detection of Edges in Spectral Data II. Nonlinear Enhancement*, SIAM Journal of Numerical Analysis, 38 (2000), 1389-1408.

[GoSh95] D. Gottlieb and C.-W. Shu, *On The Gibbs Phenomenon IV: recovering exponential accuracy in a sub-interval from a Gegenbauer partial sum of a piecewise analytic function*, Math. Comp., 64 (1995), pp. 1081-1095.

[GoSh98] D. Gottlieb and C.-W. Shu, *On the Gibbs phenomenon and its resolution*, SIAM Review 39 (1998), pp. 644-668.

[GoTa85] D. Gottlieb and E. Tadmor, *Recovering pointwise values of discontinuous data within spectral accuracy*, in "Progress and Supercomputing in Computational Fluid Dynamics", Proceedings of 1984 U.S.-Israel Workshop, Progress in Scientific Computing, Vol. 6 (E. M. Murman and S. S. Abarbanel, eds.). Birkhauser, Boston, 1985, pp. 357-375.

[HEOC85] A. Harten, B. Engquist, S. Osher and S.R. Chakravarthy, *Uniformly high order accurate essentially non-oscillatory schemes. III*, Jour. Comput. Phys. 71, 1982, 231–303.

[Jo] F. John, Partial Differential Equations,

[MMO78] A. Majda, J. McDonough and S. Osher, *The Fourier method for nonsmooth initial data*, Math. Comput. 30 (1978), pp. 1041-1081.

[Sh97] C.-W. Shu, *Essentially non oscillatory and weighted essentially non oscillatory schemes for hyperbolic conservation laws*, in "Advanced Numerical Approximation of Nonlinear Hyperbolic Equations" (A. Quarteroni, ed), Lecture Notes in Mathematics #1697, Cetraro, Italy, 1997.

[Ta94] E. Tadmor, *Spectral Methods for Hyperbolic Problems*, from "Lecture Notes Delivered at Ecole Des Ondes", January 24-28, 1994. Available at http://www.math.ucla.edu/~tadmor/pub/spectral-approximations/Tadmor.INRIA-94.pdf

[Tao] Terence Tao, Private communication.